%% file: main.tex
\documentclass[paper]{ieice}
\usepackage[dvipdfmx]{graphicx, xcolor}
\usepackage[fleqn]{amsmath}
\usepackage{amssymb, amsthm}
\usepackage{newtxtext}
\usepackage[varg]{newtxmath}
\usepackage{enumitem}
\usepackage{algorithm, algpseudocode}
\usepackage{cite}
\usepackage{color}
\usepackage{macro}

\field{A}
\title{Approximate Simultaneous Diagonalization of Matrices via Structured Low-Rank
Approximation}
\authorlist{%
  \authorentry[akema@sp.ce.titech.ac.jp]{Riku Akema}{s}{TokyoTech}\MembershipNumber{1808273}
  \authorentry[myamagi@sp.ce.titech.ac.jp]{Masao Yamagishi}{m}{TokyoTech}\MembershipNumber{0903663}
  \authorentry[isao@sp.ce.titech.ac.jp]{Isao Yamada}{f}{TokyoTech}\MembershipNumber{8510504}
}
\affiliate[TokyoTech]{The author is with the Department of Information and
Communications Engineering, Tokyo Institute of Technology 2-12-1-S3-60
Ookayama, Meguro-ku, Tokyo 152-8550, Japan.
}

\received{2020}{5}{10}
\revised{2020}{9}{7}

\begin{document}
\maketitle
\begin{summary}
  \input{Sections/abst.tex}
\end{summary}

\renewcommand\thefootnote{\arabic{footnote}}
\setcounter{footnote}{0}

\input{Sections/intr.tex}

\input{Sections/prlm.tex}

\input{Sections/atds.tex}

\input{Sections/nmex.tex}
\input{Sections/cncl.tex}

\bibliographystyle{ieicetr}
\bibliography{./BibTeXs/IEICE2020.bib}

\appendix
\input{Sections/apdx.tex}

\profile[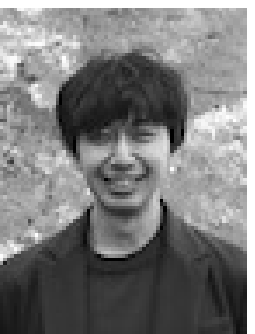]{Riku Akema}{
  received the B.E. degree in computer science and M.E. degree in information
  and communications engineering from the Tokyo Institute of Technology in 2016
  and 2018, respectively. Currently, he is a Ph.D. student in the Department of
  Information and Communications Engineering, Tokyo Institute of Technology.
  His current research interests are in signal processing and multi-way data
  analytics.
}
\profile[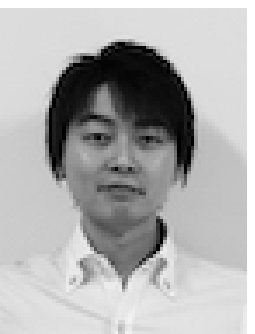]{Masao Yamagishi}{
  received the B.E. degree in computer science and the M.E. and Ph.D. degrees
  in communications and integrated systems from the Tokyo Institute of
  Technology, Tokyo, Japan, in 2007, 2008, and 2012, respectively. Currently,
  he is an assistant professor with the Department of Information and
  Communications Engineering, Tokyo Institute of Technology.
}
\profile[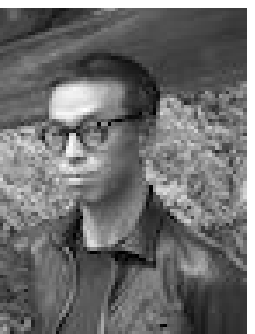]{Isao Yamada}{
  received the B.E. degree in computer science from the University
  of Tsukuba in 1985 and the M.E. and the Ph.D. degrees in electrical and
  electronic engineering from the Tokyo Institute of Technology, in 1987 and
  1990.  Currently, he is a professor with the Department of Information and
  Communications Engineering,  Tokyo Institute of Technology. He has been the
  IEICE Fellow and IEEE Fellow since 2015.
}

\end{document}

%% file: Sections/abst.tex
\emph{Approximate Simultaneous Diagonalization~(ASD)} is a problem to find a
common similarity transformation which approximately diagonalizes a given
square-matrix tuple. Many data science problems have been reduced into ASD
through ingenious modelling. For ASD, the so-called
\emph{Jacobi-like} methods have been extensively used. However, the methods
have no guarantee to suppress the magnitude of off-diagonal entries of the
transformed tuple even if the given tuple has a common exact diagonalizer,
i.e., the given tuple is \emph{simultaneously diagonalizable}.  In this paper, to
establish an alternative powerful strategy for ASD, we present a novel two-step
strategy, called \emph{Approximate-Then-Diagonalize-Simultaneously~(ATDS)}
algorithm. The ATDS algorithm decomposes ASD into (Step~1) finding a simultaneously
diagonalizable tuple near the given one; and (Step~2) finding a common similarity
transformation which diagonalizes exactly the tuple obtained in Step~1. The
proposed approach to Step~1 is realized by solving a \emph{Structured Low-Rank
Approximation~(SLRA)} with \emph{Cadzow's algorithm}. In Step~2, by
exploiting the idea in the constructive proof regarding the conditions for the
exact simultaneous diagonalizability, we obtain a common exact diagonalizer of
the obtained tuple in Step~1 as a solution for the original ASD. Unlike the
Jacobi-like methods, the ATDS algorithm has a guarantee to find a common exact
diagonalizer if the given tuple happens to be simultaneously diagonalizable.
Numerical experiments show that the ATDS algorithm achieves better performance
than the Jacobi-like methods.

%% file: Sections/intr.tex
\section{Introduction}
\label{sec:intr}
One of the most important tasks in data sciences is to extract certain common
features from given multiple data~\cite{Jolliffe2002,Bishop2006}.
Many
approaches~\cite{DeLathauwer2004,Roemer2013,Luciani2014,Andre2020,Cardoso1993,Albera2005,DeLathauwer2008,Luciani2015,VanderVeen1992,Lemma2003,Haardt2004}
reduce such tasks into a problem to find a common similarity transformation
which diagonalizes simultaneously a certain matrix tuple observed \textcolor{black}{under the influence of noise.}
A matrix tuple is called in this paper a \emph{simultaneously diagonalizable}
tuple~(see Definition~\ref{dfn:sdt}) if it has a common similarity
transformation which diagonalizes exactly all matrices in the tuple.
\textcolor{black}{
  Let $\SD$ be the set of all simultaneously diagonalizable tuples in $\Cnn
  \times \Cnn \times \cdots \times \Cnn \eqqcolon \btK \Cnn$~(see
  Definition~\ref{dfn:sdt} for the precise definition of $\SD$ and~\ref{apd:ncsd} for the nonconvexity of $\SD$). For a given $\mbA = (A_{1}, \ldots , A_{K}) \in
  \SD$, the \emph{exact simultaneous diagonalization} requires to find $S \in
  \Cnn$ such that $S^{-1} A_{k} S\ (k = 1, \ldots , K)$ are diagonal.  In this
  paper, for a general $\mbA \in \btK \Cnn$, we consider the following problem.
}
\begin{prb}[Approximate Simultaneous Diagonalization: ASD]
  \label{prb:asd}
  For a given $\mbA = (A_1, A_2, \ldots , A_K) \in \btK \Cnn$, approximate
  $\mbA$ with a certain $\tbA = (\tA_1, \ldots , \tA_K) \in \SD$ and find
  $\tS \in \Cnn$ such that $\tS^{-1} \tA_k \tS\ (k = 1, \ldots, K)$ are
  diagonal.\footnote[1]{
    In most data science
    applications~\cite{DeLathauwer2004,Roemer2013,Luciani2014,Andre2020,Cardoso1993,Albera2005,DeLathauwer2008,Luciani2015,VanderVeen1992,Lemma2003,Haardt2004}
    of ASD, the observed matrix tuple $\mbA = (A_1, \ldots , A_K) \in \btK
    \Cnn$ is modeled as a slightly perturbed version of a certain $(A^\star_1,
    \ldots , A^\star_K) \in \SD$.
  }
\end{prb}

A conventional nonconvex optimization model for the ASD of $\mbA = (A_1, \ldots
, A_K) \in \btK \Cnn$ is to find a minimizer of $f_\mbA(S) \coloneqq \sum_{k = 1}^K
\off(S^{-1} A_k S)$, where $\off \colon [x_{i, j}]_{i, j = 1}^n \mapsto \sum_{1
\le i \neq j \le n} \abs{x_{i, j}}^2$. For this optimization model, the so-called
\emph{Jacobi-like} methods have been used
extensively~\cite{Fu2006,Iferroudjene2009,Luciani2014,Luciani2015}. Although
such methods can be applied directly to general matrix tuples, they can
update only few variables in an estimate of a minimizer of $f_\mbA$ with certain
parameterized matrices, e.g., the Givens rotation matrices and the shearing
matrices. More importantly, the Jacobi-like methods have no guarantee to suppress
$f_\mbA$ even for $\mbA \in \mfD^\mathrm{abl}_{2, 2}$~(see
Example~\ref{exm:shrt} in Section~\ref{ssec:jcb} for its details). On the other
hand, the idea found in the constructive proof for a classical relation between
\emph{diagonalizability} and \emph{commutativity}~(see Fact~\ref{fct:ext})
regarding the exact simultaneous diagonalizability suggests that we can solve
the exact simultaneous diagonalization for $\mbA \in \SD$ algebraically (Note:
such an exact simultaneous diagonalization along this idea is found as,
e.g.,~\emph{Diagonalize-One-then-Diagonalize-the-Other~(DODO)} method
in~\cite{Bunse-Gerstner1993}; see Algorithm~\ref{alg:dodo} in
Section~\ref{ssec:jcb}). This situation suggests the possibility to establish a
new powerful strategy which can supersede the Jacobi-like methods for ASD if a
simultaneously diagonalizable tuple can be found as a good approximation of
$\mbA \in \btK \Cnn$.

In this paper, to establish such a computational scheme, we present a novel
two-step strategy, named
\emph{Approximate-Then-Diagonalize-Simultaneously~(ATDS)} algorithm, as a
practical solution of the following problem.
\begin{prb}[Two steps for approximate simultaneous diagonalization]
  \label{prb:2stpStrg}
  For a given $\mbA \in \btK \Cnn$,
  \begin{description}
    \item[Step 1 (Approximation):]
      approximate\footnotemark[2] $\mbA$ with a certain $\tbA = (\tA_1, \ldots , \tA_K) \in
      \SD$, where $\tbA = \mbA$ must be employed if $\mbA \in \SD$;
    \item[Step 2 (Simultaneous diagonalization):] find a common exact
      diagonalizer of $\tbA = (\tA_1, \ldots , \tA_K)$, i.e., find $\tS \in
      \Cnn$ s.t. $\tS^{-1} \tA_k \tS\ (k = 1, \ldots , K)$ are diagonal.
  \end{description}
\end{prb}
\footnotetext[2]{
  We believe that Step~1 in Problem~\ref{prb:2stpStrg} is essential for
  solving Problem~\ref{prb:asd}. However, the nonconvexity of the set
  $\SD$~(see~\ref{apd:ncsd}) makes this step intractable at least via
  straightforward approaches if we formulate Step~1 (and its relaxed
  version Step~1' in Section~3.1) as a certain nonconvex feasibility
  problem or nonconvex optimization problem under a certain specified
  criterion.  Indeed, we have not yet found such guaranteed algorithms
  achieving $\tbA \in \SD$ within a certain ball, of a prescribed radius,
  centered at $\mbA$ and thus propose in this paper to relax Step~1 in
  Problem~\ref{prb:2stpStrg} by Step~1' with \emph{Structured Low-Rank
  Approximation~(SLRA)}~(see Section~3.1).  To promote further
  breakthrough for various innovative approximations,  we formulate
  Step~1 in Problem~\ref{alg:ovralg} (and its relaxed version Step~1' in
  Section~3.1) without specifying approximation criteria.
}

Since Step~2 can be solved by the DODO method, we focus on Step~1. The proposed
approach to Step~1 is designed based on the fact that a necessary condition for
$\mbX \in \SD$ can be translated into a certain rank condition for a structured
matrix defined with $\mbX$~(see Theorem~\ref{thr:XiAndSD}). More precisely, 
\textcolor{black}{
  we propose to relax Step~1 by Step~1' with \emph{Structured Low-Rank
  Approximation~(SLRA)}~\cite{Markovsky2008} and then propose to solve Step~1'
  with \emph{Cadzow's algorithm}~\cite{Cadzow1988}. Unlike the Jacobi-like
  methods, the proposed ATDS algorithm has a guarantee to find a common exact
  diagonalizer if $\mbA$ happens to satisfy $\mbA \in \SD$. Numerical
  experiments show that, compared with the Jacobi-like methods, the proposed
  ATDS algorithm achieves a better approximation to the desired common
  similarity transformation at the expense of reasonable computational time.
}

\textcolor{black}{
  A preliminary version of this paper was presented at a conference~\cite{Akema2019}.
}

\noindent {\bf Notations (see also Table~\ref{tbl:ntn}).} 
\textcolor{black}{
  Let $\N$, $\R$, and $\C$ denote the set of all nonnegative integers, the set
  of all real numbers, and the set of all complex numbers, respectively.
}
For $\mbx \in \C^n$, $\|\mbx\|$ denotes the Euclidean norm of $\mbx$. The
matrix $I_n$ denotes the $n$-by-$n$ identity matrix. For $X \in \Cmn$,
$X^\top$, $X^\ast$, $X^{-1}$, $X^{-\top}$, $\trc(X)$, $\|X\|_F$, \textcolor{black}{$\range(X)$}, and $\nullspace(X)$
denote \textcolor{black}{respectively the transpose, the conjugate transpose, the inverse, the inverse of transpose, the trace,
the Frobenius norm, the range space, and the nullspace of $X$.} The mapping $\vec :
\Cmn \rightarrow \C^{m n}$ denotes vectorization by stacking the columns of a
matrix and $\vec^{-1}$ its inverse. For $\mbx \in \Cn$, $\diag(\mbx) \in \Cnn$
is a diagonal matrix whose diagonal entries are given by the components of
$\mbx$.  The Kronecker product of $X = [x_{i, j}] \in \Cmn$ and $Y \in \C^{p
\times q}$ is $X \otimes Y \in \C^{m p \times n q}$ whose $(i, j)$th blocks are
$[X \otimes Y]_{i, j} = x_{i,j} Y \ (i = 1, 2, \ldots , m; j = 1, 2, \ldots ,
n)$.  The Khatri-Rao product of $X = [\mbx_1 \ \cdots \ \mbx_l] \in \C^{m
\times l}$ and $Y = [\mby_1 \ \cdots \ \mby_l] \in \C^{n \times l}$ is $X \odot
Y \coloneqq [\mbx_1 \otimes \mby_1 \ \cdots \ \mbx_l \otimes \mby_l] \in
\C^{mn\times l}$. \textcolor{black}{The direct sum of
$\check{X}_p \in \C^{n_p \times n_p} \ (p = 1, \ldots , d; n_1 + \cdots + n_d =
n)$ is the block diagonal matrix $\check{X}_1 \oplus \check{X}_2 \oplus \cdots \oplus \check{X}_d \in
\Cnn$ with the $i$th diagonal block $\check{X}_i~(i = 1, \ldots , d)$.}  For $X_p \in \C^{m \times n} \ (p = 1,
\ldots , d)$, $\spn\{X_1, \ldots , X_d\} \subset \C^{m \times n}$ is the set of
all linear combinations of $X_1, \ldots , X_d$. 
\textcolor{black}{
  For $\mbX = (X_{1}, \ldots , X_{K}) \in \btK \Cnn$, we define
  $\norm{\mbX}_{F} \coloneqq \sqrt{\sum_{k=1}^{K}\norm{X_{k}}_{F}^{2}}$.
}

%% file: Sections/prlm.tex
\section{Preliminaries}
\label{sec:prlm}
\subsection{\textcolor{black}{Useful Facts on Diagonalizability}}
\label{ssec:dfnfct}
Recall that $X \in \Cnn$ is said to be \emph{diagonalizable} if there
exists $S \in \Cnn$ such that $S^{-1} X S$ is diagonal. In the following, we use
$\Dga \subset \Cnn$ to denote the set of all diagonalizable matrices.
The simultaneous diagonalizability in Definition~\ref{dfn:sdt} below is a
natural extension of the diagonalizability of a matrix \textcolor{black}{and defined} for a matrix tuple.
\begin{dfn}[Simultaneously diagonalizable tuple~{\cite[Definition~1.3.20]{Horn2013}}]
  \label{dfn:sdt}
  A matrix tuple $\mbX = (X_1, \ldots , X_K) \in \btK \Cnn$ is said to be
  \emph{simultaneously diagonalizable} if there exists a common $S \in \Cnn$
  s.t. $S^{-1} X_k S\ (k = 1, \ldots , K)$ are diagonal. In this paper, we use
  $\SD$ to denote the set of all simultaneously diagonalizable tuples in $\btK
  \Cnn$.
\end{dfn}

The exact simultaneous diagonalization of $\mbX \in \SD$ is said to be
\emph{essentially unique}~\cite{DeLathauwer2004} if its exact common
diagonalizer $S \in \Cnn$ is determined uniquely up to a permutation and a
scaling of the column vectors of $S$. The following is well-known as an
equivalent condition for the essential uniqueness.
\begin{fct}[Neccesary and sufficient condition for essential uniqueness{~\cite[Theorem~6.1]{DeLathauwer2004}}]
  \label{fct:essUnq}
  \textcolor{black}{Suppose that $S \in
    \Cnn$ is a common diagonalizer of $\mbX = (X_{1}, \ldots , X_{K}) \in \SD$ satisfying $\Lambda_{k} \coloneqq S^{-1} X_k S =
    \diag(\lambda_1^{(k)}, \ldots , \lambda_n^{(k)})$, where
    $\lambda^{(k)}_{1}, \ldots, \lambda^{(k)}_{n} \in \C$ are all
  eigenvalues of $X_{k}\ (k = 1, \ldots , K)$ (Note: The order of eigenvalues $\lambda^{(k)}_{1}, \ldots , \lambda^{(k)}_{n}$ is determined according to $S$).} Then, the exact simultaneous
  diagonalization of $\mbX$ is essentially unique if and only if
  $[\lambda_p^{(1)}\ \lambda_p^{(2)}\ \cdots\ \lambda_p^{(K)}] \neq
  [\lambda_q^{(1)}\ \lambda_q^{(2)}\ \cdots\ \lambda_q^{(K)}]$ for any $p \neq
  q; p, q \in \{1, \ldots , n\}$.
\end{fct}

As seen below, the simultaneous diagonalizability of a matrix tuple requires
not merely the diagonalizability of every matrix but also the commutativity of
all matrices therein.
\begin{fct}[{Necessary and sufficient condition for $\SD$~\cite[Theorem 1.3.21]{Horn2013}}]
  \label{fct:ext}
  Let $\Cmm \subset \btK \Cnn$ denote the set of all tuples $(X_1, \ldots
  , X_K)$ such that $X_k X_l = X_l X_k$, i.e., $X_k$ and $X_l$ commute, for any
  $k, l \in \{1, \ldots , K\}$. Then,
  \begin{align}
    \label{eq:nssd}
    \SD = \Cmm \cap (\Dga)^K.
  \end{align}
  That is, a tuple is (exactly) simultaneously diagonalizable if and only if every pair
  $X_k$ and $X_l$ commutes and every $X_k$ is diagonalizable.  
\end{fct}
\begin{rmk}
  The constructive proof for Fact~\ref{fct:ext} found in, e.g.,~\cite{Horn2013}
  can be translated into a finite-step algorithm (the DODO
  method~\cite{Bunse-Gerstner1993}) for Step~2 in
  Problem~\ref{prb:2stpStrg}~(see Section~\ref{ssec:dodo}).
\end{rmk}

\begin{table}[t]
  \centering
  \caption{Notations on the sets used in this paper.}
  \label{tbl:ntn}
  \begin{tabular}{ll}
    \hline
    \textcolor{black}{$\N$} & \textcolor{black}{the nonnegative integers} \\
    $\R$ & the real numbers \\
    $\R^{m \times n}$ & $m$-by-$n$ real matrices \\
    $\C$ & the complex numbers \\
    $\C^{m \times n}$ & $m$-by-$n$ complex matrices \\
    $\Dga$ & diagonalizable matrices in $\Cnn$ \\ 
    $\Cmm$ & commuting tuples in $\btK \Cnn$ \\ 
    $\SD$ & simultaneously diagonalizable tuples in $\btK \Cnn$ \\ 
    \hline
  \end{tabular}
\end{table}

\subsection{Jacobi-like methods}
\label{ssec:jcb}
To see the basic idea of the Jacobi-like methods, let us explain the scheme of
sh-rt~\cite{Fu2006}, for $\mbA = (A_1, \ldots , A_K) \in \btK \Rnn$, which is
known as the one of the earliest extensions, of the Jacobi
methods~\cite{Jacobi1846}~(see also, e.g.,~\cite[Section~8.4]{Golub1996}) for a
symmetric-matrix diagonalization, to simultaneous diagonalization. Sh-rt uses
the shearing matrix $H(p, q, \phi) = [h_{i, j}] \in \Rnn$ and the Givens
rotation matrix $G(p, q, \theta) = [g_{i, j}] \in \Rnn\ (p, q \in \{1, \ldots ,
n\}; p \neq q)$ whose entries are the same with $I_n$ except for $h_{p, p} =
h_{q, q} = \cosh \phi; h_{p, q} = h_{q, p} = \sinh \phi$ and $g_{p, p} = g_{q,
q} = \cos \theta; g_{p, q} = -g_{q, p} = -\sin \theta$, respectively.

Set \textcolor{black}{$\breve{S} \coloneqq I_{n}$ and} $\breve{A}_k = [\breve{a}^{(k)}_{i, j}] \coloneqq A_k\ (k = 1, \ldots ,
K)$.  Then sh-rt repeats the following procedures: (i)~choose $p$ and $q\ (p < q)$; (ii)~find
$\phi^\star$ to enhance normality of $(H(p, q, \phi^\star))^{-1} \breve{A}_k
H(p, q, \phi^\star) \eqqcolon [\breve{b}^{(k)}_{i, j}] = \breve{B}_k \in \Rnn\
(k = 1, \ldots , K)$ with $\breve{\mbB} \coloneqq (\breve{B}_1, \ldots,
\breve{B}_K)$; (iii)~find $\theta^\star$ for suppression of
$\mcG_{\breve{\mbB}}(\theta) \coloneqq \sum_{k = 1}^n \off((G(p, q,
\theta))^{-1} \breve{B}_k G(p, q, \theta))$; and \textcolor{black}{(iv)~set $\breve{A}_k \coloneqq
(G(p, q, \theta^\star))^{-1}\breve{B}_k G(p, q, \theta^\star), \breve{\Lambda}_{k} \coloneqq \diag(\bra^{(k)}_{1,1}, \ldots , \bra^{(k)}_{n,n})\ (k = 1, \ldots , K)$ and $\breve{S} \coloneqq \breve{S} H(p, q, \phi^\star) G(p, q, \theta^{\star})$.}

After lengthy algebra, the authors of~\cite{Fu2006} suggest to use
$\phi^\star$ and $\theta^\star$ satisfying the following conditions with $l \in
\argmax_{k \in \{1, \ldots , K\}} \abs{\breve{a}^{(k)}_{p, p} -
\breve{a}^{(k)}_{q, q}}$:
\begin{align}
  \label{eq:shrt1}
  & \tanh \phi^\star = \frac{\kappa^{(l)}_{p, q} - (\breve{a}^{(l)}_{p, p} - \breve{a}^{(l)}_{q, q}) (\breve{a}^{(l)}_{p, q} - \breve{a}^{(l)}_{q, p})}{2 ((\breve{a}^{(l)}_{p, p} - \breve{a}^{(l)}_{q, q})^2 + (\breve{a}^{(l)}_{p, q} - \breve{a}^{(l)}_{q, p})^2) + \xi^{(l)}_{p, q}} \\
  \label{eq:shrt2}
  & \tan 4 \theta^\star = \frac{2 \sum_{k = 1}^K (\breve{b}^{(k)}_{p, p} - \breve{b}^{(k)}_{q, q}) (\breve{b}^{(k)}_{p, q} + \breve{b}^{(k)}_{q, p})}{\sum_{k = 1}^K (\breve{b}^{(k)}_{p, p} - \breve{b}^{(k)}_{q, q})^2 - (\breve{b}^{(k)}_{p, q} + \breve{b}^{(k)}_{q ,p})^2},
\end{align}
where $\kappa^{(l)}_{p, q} \coloneqq \sum_{j = 1; j \neq p, q}^n
(\breve{a}^{(l)}_{p, j} \breve{a}^{(l)}_{q, j} - \breve{a}^{(l)}_{j, p}
\breve{a}^{(l)}_{j, q})$ and $\xi^{(l)}_{p, q} \coloneqq \sum_{j = 1; j \neq p,
q}^n (\breve{a}^{(l)}_{p, j})^2 + (\breve{a}^{(l)}_{q, j})^2 +
(\breve{a}^{(l)}_{j, p})^2 + (\breve{a}^{(l)}_{j, q})^2$.  Indeed,
$\theta^\star$ in~\eqref{eq:shrt2} is a stationary point of $\mcG_{\bbB}$,
which is verified in~\cite{Fu2006,Iferroudjene2009}. However, as seen in
Example~\ref{exm:shrt} below, we remark that sh-rt has no guarantee to suppress
$f_\mbA(S) = \sum_{k = 1}^K \off(S^{-1} A_k S)$ even for $\mbA \in
\mfD^\mathrm{abl}_{2, 2}$.
\begin{exm}[A weakness of a Jacobi-like method: sh-rt]
  \label{exm:shrt}
  Suppose that $A_1 \coloneqq [a_{i, j}] \in \R^{2 \times 2}$ and $A_2 = c_0
  I_2 + c_1 A_1 \in \R^{2 \times 2}$, where $a_{1, 1} = a_{2, 2},\ \abs{a_{1,
  2}} \neq \abs{a_{2, 1}},\ a_{1, 2} a_{2, 1} > 0$, and $c_0, c_1 \neq 0$.
  Since the discriminant of the characteristic polynomial of $A_1$ is positive,
  $A_1$, and then $A_2$, have distinct real eigenvalues and therefore they are
  diagonalizable.  By commutativity of $A_1$ and $A_2$, moreover, we see $(A_1,
  A_2) \in \mfD^\mathrm{abl}_{2, 2}$~(see Fact~\ref{fct:ext}).

  Apply sh-rt scheme to $(A_1, A_2)$ with $(p, q) = (1, 2)$ and $l = 1$. Then,
  since $\kappa^{(1)}_{1, 2} = 0$ and $\breve{a}^{(1)}_{1, 1} -
  \breve{a}^{(1)}_{2, 2} = 0$ in~\eqref{eq:shrt1}, we get $\tanh \phi^\star =
  0$, i.e., $H(1, 2, \phi^\star) = I_2$. Moreover, since $\breve{b}^{(k)}_{1,
  1} - \breve{b}^{(k)}_{2, 2} = \breve{a}^{(k)}_{1, 1} - \breve{a}^{(k)}_{2, 2}
  = 0\ (k = 1, 2)$ in~\eqref{eq:shrt2}, we get $\tan 4 \theta^\star = 0$, i.e.,
  $G(1, 2, \theta^\star) = I_2$. Therefore, since $\breve{A}_k = \breve{B}_k =
  A_k\ (k = 1, 2)$, we see $f_{(\breve{A}_1, \breve{A}_2)} = f_{(A_1, A_2)}$.
\end{exm}

\subsection{The DODO Method}
\label{ssec:dodo}
\begin{algorithm}[t]
  \caption{The DODO method~\cite{Bunse-Gerstner1993}}
  \label{alg:dodo}
  \begin{algorithmic}
    \Function{DODO}{$(A_1, \ldots , A_K), n$}
    \If{$A_k\ (k = 1, \ldots , K)$ are diagonal}
    \State \Return $I_n$
    \Else
    \State Diagonalize $A_l\ (l \in \{1, \ldots , K\})$ as $S_0^{-1} A_l S_0$ of the form~\eqref{eq:dsDg}.
    \State $S_0^{-1} A_k S_0 = \check{A}^{(k)}_1 \oplus \cdots \oplus \check{A}^{(k)}_d\ (k = 1, \ldots , K)$
    \For{$p = 1$ to $d$}
    \State $S_p =$ \Call{DODO}{$(\check{A}^{(1)}_p, \ldots , \check{A}^{(K)}_p), n_p$}
    \EndFor
    \State \Return $S = S_0 (S_1 \oplus \cdots \oplus S_d)$
    \EndIf
    \EndFunction
  \end{algorithmic}
\end{algorithm}

The key idea behind the DODO method (Algorithm~\ref{alg:dodo}) is to reduce the
exact simultaneous diagonalization of $\mbA$ into exact simultaneous
diagonalizations of tuples of smaller matrices.  To see how Algorithm~\ref{alg:dodo}
works, let us demonstrate its procedures.
\begin{enumerate}[label=(\roman*)]
  \item Choose $l \in \{1, \ldots , K\}$ s.t. $A_l \in \Dga$ is not diagonal,
    arbitrarily. Diagonalize $A_l$ with $S_0 \in \Cnn$ as
    \begin{align}
      \label{eq:dsDg}
      S_0^{-1} A_l S_0 = \lambda_1^{(l)} I_{n_1} \oplus \cdots \oplus \lambda_d^{(l)} I_{n_d},
    \end{align}
    where $\lambda_1^{(l)}, \ldots , \lambda_d^{(l)} \in \C$ are
    distinct\footnote[3]{Such $S_0$ can be computed by permuting column vectors of
    any $S$ satisfying $S^{-1} A_l S$ is diagonal.} and $\max\{n_1, \ldots ,
    n_d\} < n$ with $n_1 + \cdots + n_d = n$ holds because $A_l$ is not a
    constant multiple of $I_n$ (Note: if $n_1 = \cdots = n_d = 1$, i.e., $d = n$
    holds, $S_0$ is a common exact diagonalizer for $\mbA$, which is verified
    essentially in the following procedure~(ii)).
  \item Compute $S_0^{-1} A_k S_0 \ (k = 1, \ldots , K; k \neq l)$ which can be expressed as
    $S_0^{-1} A_k S_0 = \check{A}^{(k)}_1 \oplus \cdots \oplus
    \check{A}^{(k)}_d$ with some $\check{A}^{(k)}_p \in
    \C^{n_p \times n_p}$ for any $p = 1, \ldots , d$ (Note: this is verified by applying
    Fact~\ref{fct:cmm}(c) to the commutativity\footnote[4]{The commutativity can be verified
      by $(S_0^{-1} A_k S_0) (S_0^{-1} A_l S_0) = S_0^{-1} A_k A_l S_0 =
      S_0^{-1} A_l A_k S_0 = (S_0^{-1} A_l S_0) (S_0^{-1} A_k S_0)$~(see
      Fact~\ref{fct:ext}).
    } of $S_0^{-1} A_k S_0$ and $S_0^{-1} A_l S_0$ in~\eqref{eq:dsDg}).
  \item Construct $\check{\mbA}_p \coloneqq (\check{A}^{(1)}_p, \ldots ,
    \check{A}^{(K)}_p) \in \btK \C^{n_p \times n_p}$ with $\check{A}^{(l)}_p
    \coloneqq \lambda_p^{(l)} I_{n_p}\ (p = 1, \ldots , d)$. Then, we
    obtain $\check{\mbA}_p \in
    \mfD^{\mathrm{abl}}_{n_p, K}$ for all $p = 1, \ldots , d$ (Note on
    \emph{Commutativity}:~for any pair $(k_1, k_2) \in \{1, \ldots , K\}^2$,
    the commutativity of $S_0^{-1} A_{k_1} S_0 = \check{A}^{(k_1)}_1 \oplus
    \cdots \oplus \check{A}^{(k_1)}_d$ and $S_0^{-1} A_{k_2} S_0 =
    \check{A}^{(k_2)}_1 \oplus \cdots \oplus \check{A}^{(k_2)}_d$ (see~(ii) for these structures)
    implies commutativity of $\check{A}^{(k_1)}_p$ and $\check{A}^{(k_2)}_p\ (p
    = 1, \ldots , K)$; on \emph{Diagonalizability}: the diagonalizability of
    $S_0^{-1} A_k S_0$~(guaranteed by Fact~\ref{fct:ext}) implies diagonalizability of
    $\check{A}^{(k)}_p\ (p = 1, \ldots , d)$~(see Fact~\ref{fct:bldg})).
  \item For each $\check{\mbA}_p\ (p = 1,
    \ldots , K)$, repeat~(i-iii), where many $K$-tuples of smaller matrices may appear in the process, until all $\check{A}^{(k)}_p\ (k = 1, \ldots ,
    K; p = 1, \ldots , d)$ become diagonal.
  \item Compute a common exact diagonalizer for $\mbA$ as $S \coloneqq S_0 (S_1
    \oplus \cdots \oplus S_d) \in \Cnn$, where $S_p \in \C^{n_p \times n_p}$ is
    a common exact diagonalizer of $\check{\mbA}_p$, for each $p = 1, \ldots ,
    d$, constructed with a diagonalizer in~(i) regarding $\check{\mbA}_p$.
\end{enumerate}

%% file: Sections/atds.tex
\section{Approximate-Then-Diagonalize-Simultaneously Algorithm}
\label{sec:trns}
\subsection{Simultaneous Diagonalizability Condition in terms of the Kronecker Sums}
It is not hard to see that $X \in \Cnn$ and $Y \in \Cnn$ commute if and only if
$\vec(Y) \in \nullspace(I_n \otimes X - X^\top \otimes I_n)$, where $I_n \otimes X -
X^\top \otimes I_n$ is called the \emph{Kronecker sum} of $X$ and $-X^\top$.
This simple fact motivates us to introduce a linear mapping $\Xi \colon \btK
\Cnn \to \C^{Kn^2\times n^2}, \mbX = (X_1, X_2, \ldots , X_K) \mapsto$
\begin{align}
  \label{eq:Xi}
  \Xi(\mbX) \coloneqq
  \left[
    \begin{array}{c}
      I_n \otimes X_1-X_1^\top\otimes I_n\\
      I_n \otimes X_2-X_2^\top\otimes I_n\\
      \vdots\\
      I_n \otimes X_K - X_K^\top \otimes I_n \\
    \end{array}
  \right] \in \C^{K n^2\times n^2}.
\end{align}
Moreover, for $\hX \in \Xi(\btK \Cnn) \coloneqq \{\Xi(\mbY) \in \C^{K n^2
\times n^2} \mid \mbY \in \btK \Cnn\}$, we introduce an affine subspace
$\Xi^{-1}(\hX) \coloneqq \{\mbY \in \btK \Cnn \mid \Xi(\mbY) = \hX\} \subset
\btK \Cnn$.

\begin{thr}[Characterizations of $\Cmm$ and $\SD$ with $\Xi$]
  \label{thr:XiAndSD}
  Let $\mbX = (X_1, \ldots ,X_K) \in \btK \Cnn$.
  \begin{enumerate}[label=(\alph*)]
    \item $\mbX \in \Cmm \Leftrightarrow (\forall k \in \{1, \ldots , K\})\
      \Xi(\mbX) \vec(X_k) = \zeros$.
    \item $\mbX \in \SD \Leftrightarrow (\exists S \in \Cnn)\ \nullspace(\Xi(\mbX))
      \supset \range(S^{-\top} \odot S)$.
    \item $\mbX \in \SD \Rightarrow
      \left\{
        \begin{array}{l}
          \rank(\Xi(\mbX))\le n^2-n, \\
          \rank(\Xi(\mbX))= n^2-n \Leftrightarrow \\
          \quad \text{exact simultaneous
        diagonalization} \\
        \qquad \text{of }\mbX \text{ is essentially unique}.
        \end{array}
      \right.$
    \item If at least one $X_l$ has $n$ distinct eigenvalues, \\
      $\mbX\in\SD\Leftrightarrow\rank(\Xi(\mbX))=n^2-n$.
    \item Let $\mfL_{n^2 - n} \coloneqq \{\hX \in \C^{K n^2
      \times n^2} \mid \rank(\hX) \le n^2 - n\}$. Suppose $\hX \in \Xi(\btK
      \Cnn) \cap \mfL_{n^2 - n}$ and $\Xi(\mbY) = \hX$ for
      some $\mbY \coloneqq (Y_1, \ldots , Y_K) \in \btK \Cnn$, where at least
      one $Y_l$ has $n$ distinct eigenvalues and is diagonalizable as $S^{-1}
      Y_l S$ with $S \in \Cnn$. Then, $\mbY \in \SD$ and $S^{-1} Y_k S\ (k = 1,
      \ldots , K)$ are diagonal.
  \end{enumerate}
\end{thr}
(The proof of Theorem~\ref{thr:XiAndSD} is given in~\ref{prf:XiAndSD}.)

\begin{rmk}[On Theorem~\ref{thr:XiAndSD}]
  \ 
  \begin{enumerate}[label=(\alph*)]
    \item Theorem~\ref{thr:XiAndSD}(c) implies $\Xi^{-1}(\Xi(\btK \Cnn)
      \cap \mfL_{n^2 - n}) \coloneqq \{\mbX \in \btK \Cnn \mid \rank(\Xi(\mbX))
      \le n^2 - n\} \supset \SD$.
    \item Theorem~\ref{thr:XiAndSD}(e) is a sufficient condition for $\mbX \in
      \SD$ in terms of the range of $\Xi$. Although this condition is similar
      to~(d), the condition in~(e) directly motivates us to solve a Structured
      Low-Rank Approximation~(SLRA) below (see Problem~\ref{prb:asdslra}) for
      finding a simultaneously diagonalizable tuple near $\mbA$.
  \end{enumerate}
\end{rmk}

For $\hX \in \Xi(\btK \Cnn)$, the projection onto $\Xi^{-1}(\hX)$ can be
computed as in Proposition~\ref{prp:pjXi}.
\begin{prp}[On projection onto $\Xi^{-1}(\hX)$]
  \label{prp:pjXi}
  Let $\hX\in\Xi(\btK \Cnn)$. Choose $\mbX^\diamond \coloneqq (X_1^\diamond,
  \ldots ,X^\diamond_K)\in\Xi^{-1}(\hX)$ arbitrarily. Then,
  \begin{enumerate}[label=(\alph*), leftmargin=*]
    \item $\Xi^{-1}(\hX)=\mbX^\diamond + (\spn\{I_n\})^K$;
    \item the projection of $\mbX = (X_1, \ldots , X_K) \in \btK \Cnn$ onto
      $\Xi^{-1}(\hX)$, i.e., $P_{\Xi^{-1}(\hX)} \colon (X_1, \ldots ,
      X_K)\mapsto$
      \begin{align*}
        (Z_1, \ldots , Z_K) \coloneqq \argmin_{(Y_1,\ldots ,Y_K)\in\Xi^{-1}(\hX)}\sum_{k=1}^K \norm{X_k - Y_k}_F^2,
      \end{align*}
      is given by $Z_k = X^\diamond_k + (\trc(X_k - X_k^\diamond) / n) I_n\
      (k = 1, \ldots , K)$;
    \item for $\mbX = (X_1, \ldots , X_K)$ and $P_{\Xi^{-1}(\hX)}(\mbX) = (Z_1,
      \ldots , Z_K)$, $\sum_{k = 1}^K \norm{X_k - Z_k}^2_F = \norm{\Xi(\mbX) -
      \hX}^2_F / 2n$. 
  \end{enumerate}
\end{prp}
(The proof of Proposition~\ref{prp:pjXi} is given in~\ref{prf:pjXi}.)

Finally, by using Theorem~\ref{thr:XiAndSD}(e), Proposition~\ref{prp:pjXi}(b),
and~\ref{prp:pjXi}(c), we propose to relax Step~1 in Problem~\ref{prb:2stpStrg}
by the following Step~1'~(see Remark~2(a)).
\begin{description}
  \item[Step~1' (Approximation with a structured low-rankness):] \textcolor{black}{approximate
    $\mbA$ with a certain $\tbA \in \Xi^{-1}(\Xi(\btK \Cnn) \cap \mfL_{n^2 -
  n})$, where $\tbA = \mbA$ must be employed if $\mbA \in \SD$;}
\end{description}

Proposition~1 suggests that Step~1' can be decomposed further into the
following Step~1'a and Step~1'b. 
\begin{prb}[Two steps for ATDS algorithm with a Structured Low-Rank Approximation: SLRA]
  \label{prb:asdslra}
  For a given $\mbA \in \btK \Cnn$,
  \begin{description}
    \item[Step 1'a (A structured low-rank approximation):]\ \\
      \textcolor{black}{approximate $\Xi(\mbA)$ with a certain $\hA_\ast \in \Xi(\btK \Cnn) \cap
      \mfL_{n^2 - n}$, where $\hA_\ast = \Xi(\mbA)$ must be
      employed if $\Xi(\mbA) \in \Xi(\btK \Cnn) \cap \mfL_{n^2 - n}$;}
    \item[Step 1'b (Projection onto an affine subspace):]\ \\ compute
      $\tbA=P_{\Xi^{-1}(\hA_\ast)}(\mbA)$~(see Proposition~\ref{prp:pjXi});
    \item[Step~2 (Simultaneous diagonalization):] find a common exact
      diagonalizer of $\tbA = (\tA_1, \ldots , \tA_K)$, i.e., find $\tS \in
      \Cnn$ s.t. $\tS^{-1} \tA_k \tS\ (k = 1, \ldots , K)$ are diagonal.
  \end{description}
\end{prb}
\begin{rmk}[On the proposed approach]
  \ 
  \begin{enumerate}[label=(\alph*)]
    \item Proposition~\ref{prp:pjXi}(c) ensures that $\tbA = (\tA_1,\ldots
      ,\tA_K) := P_{\Xi^{-1}(\hA_\ast)}(\mbA)$ is expected to be close to
      $\mbA$ \textcolor{black}{if $\hA_{\ast}$ is obtained in Step~1'a as a good approximation of
      $\Xi(\mbA)$.} From Theorem~\ref{thr:XiAndSD}(e), moreover, $\tbA$ is
      guaranteed to be simultaneously diagonalizable if $\tA_l$ has $n$
      distinct eigenvalues for some $l \in \{1, \ldots , K\}$.
    \item The proposed approach exploits an algebraic property, i.e., $\SD =
      \Cmm \cap (\Dga)^K$, of a simultaneously diagonalizable tuple. Using this
      algebraic property aims to achieve a denoising effect in Step~1'. The
      effectiveness, of using the commutativity condition, for denoising in ASD
      was suggested in~\cite{VanderVeen1992} but only for $K = 2$.
  \end{enumerate}
\end{rmk}

\subsection{Approximate Simultaneous Diagonalization Algorithm by Cadzow's Algorithm}
\label{ssec:palg}
We have already shown how to solve Step 1'b in
Proposition~\ref{prp:pjXi}(b).\footnote[5]{
  For $\hX \in \Xi(\btK \Cnn)$, although the choice of $\mbX^\diamond \in
  \Xi^{-1}(\hX)$ in Proposition~\ref{prp:pjXi} is arbitrary, a possible and
  simple one is the following.  Partition $\hX = [\hX_1^\top\ \cdots \
  \hX_K^\top]^\top$, where $\hX_k \coloneqq I_n \otimes Y_k - Y_k^\top \otimes I_n \in
  \C^{n^2 \times n^2}\ (k = 1, \ldots , K)$ for some $Y_k = [y^{(k)}_{i, j}]
  \in \Cnn$. Since the top-left blocks of $\hX_k\ (k = 1, \ldots , K)$ are $Y_k
  - y^{(k)}_{1, 1} I_n$ and then satisfy $I_n \otimes (Y_k - y^{(k)}_{1, 1}
  I_n) - (Y_k - y^{(k)}_{1, 1} I_n)^\top \otimes I_n = \hX_k$, we can use $(Y_1
  - y_{1, 1}^{(1)} I_n, \ldots , Y_K - y^{(K)}_{1, 1} I_n)$ as $\mbX^\diamond$.
} To realize Step 1'a, we propose to use Cadzow's
algorithm~\cite{Cadzow1988} also known as alternating projection algorithm,
below:
\begin{align}
  \label{eq:pral}
  \left\{
    \begin{array}{l}
      \hA(0) \coloneqq \Xi(\mbA); \\
      \hA(t + 1) \coloneqq P_{\Xi(\btK \Cnn)} \circ P_{\mfL_{n^2 - n}}(\hA(t))
    \end{array}
  \right.
\end{align}
$(t = 0, 1, \ldots)$, where $\circ$ denotes a composition of mappings and, for
any $\hX \in \C^{K n^2 \times n^2}$,
\begin{align}
  \label{eq:prjs}
  \left\{
    \begin{aligned}
      & P_{\mfL_{n^2 - n}}(\hX) \in \argmin_{\hY \in \mfL_{n^2 - n}} \norm{\hX - \hY}_F; \\
      & P_{\Xi(\btK \Cnn)}(\hX) \coloneqq \argmin_{\hY \in \Xi(\btK \Cnn)} \norm{\hX - \hY}_F.
    \end{aligned}
  \right.
\end{align}

\begin{prp}[Monotone approximation property of alternating projection]
  \label{prp:mnap}
  Let $(\hA(t))_{t = 0}^\infty$ be the sequence generated by~\eqref{eq:pral}.
  Then, the sequence $(\tbA(t))_{t=0}^\infty$ defined by
  $\tbA(t):=P_{\Xi^{-1}(\hA(t))}(\mbA)\in\btK \Cnn$ satisfies
  $\Xi(\tbA(t))=\hA(t)\in\Xi(\btK \Cnn)\ (t=0,1,\ldots)$ and
  \begin{align*}
    &\norm{\Xi(\tbA(t+1))-P_{\mfL_{n^2-n}}(\Xi(\tbA(t+1)))}_F\\
    &\le\norm{\Xi(\tbA)(t))-P_{\mfL_{n^2-n}}(\Xi(\tbA(t)))}_F\quad(t=0,1,\ldots).
  \end{align*}
\end{prp}
(The proof of Proposition~\ref{prp:mnap} is given in~\ref{prf:mnap}.)

Finally, we propose the ATDS algorithm with SLRA (Algorithm~\ref{alg:ovralg}),
as a practical solution to Problem~\ref{prb:asdslra}, 
  where we also propose additionally to use "Pseudo Common Diagonalizer~(PCD)"~(see function
  \texttt{PCD} in Algorithm~\ref{alg:ovralg}) under Assumption~1\footnotemark[6] for $P_{\Xi^{-1}(\hA(t_{\mathrm{end}}))}(\mbA)$~(see Remark~\ref{rmk:prpalg}(e)) if $\tbA \in \SD$ is not achieved by Step~1'b.
  \footnotetext[6]{
  \textcolor{black}{
    Assumption~1 for $P_{\Xi^{-1}(\hA(t_{\mathrm{end}}))}(\mbA)$ seems to be weak
    enough in practice as we have not seen any exceptional case in our numerical
    experiments~(see Section~\ref{sec:nuex}).
  }
}
\textcolor{black}{
  Note that, if $\mbA \in \SD$,
Algorithm~\ref{alg:ovralg} has guarantee to satisfy $\tbA = \mbA$~(a requirement in Step~1 of Problem~2) and to find the common exact diagonalizer
of $\mbA$ with the DODO method~(see Algorithm~\ref{alg:dodo}) unlike the Jacobi-like methods~(see Example~1 in Section~\ref{ssec:jcb}).
}
\begin{prp}[On Algorithm~\ref{alg:ovralg} in the case of $\mbA \in \SD$]
  \label{prp:slvesd}
  Suppose that $\mbA = (A_1, \ldots , A_K) \in \btK \Cnn$ happens to satisfy
  $\mbA \in \SD$. Then, Algorithm~\ref{alg:ovralg} finds a common exact
  diagonalizer $\tS \in \Cnn$ s.t.  $\tS^{-1} A_k \tS\ (k = 1, \ldots , K)$ are
  diagonal.
\end{prp}
\begin{proof}
  Recall that an initial guess $\hA(0) = \Xi(\mbA)$ with $\mbA \in \SD$
  satisfies $\hA(0) \in \Xi(\btK \Cnn) \cap \mfL_{n^2 - n}$~(see Remark~2(a)),
  which implies that $\norm{\hA(0) - P_{\mfL_{n^2 - n}}(\hA(0))}_F = \norm{O}_F
  = 0$. Therefore, we get $\hA_\ast \coloneqq \hA(0) = \Xi(\mbA)$ by Cadzow's step and then
  $\tbA \coloneqq P_{\Xi^{-1}(\hA_\ast)}(\mbA) = \mbA \in \SD$. Consequently, $\tS \in \Cnn$
  found by the DODO method makes $\tS^{-1} A_k \tS\ (k = 1,
  \ldots , K)$ diagonal.
\end{proof}

\newcommand{\algrule}[1][.2pt]{\par\vskip.5\baselineskip\hrule height #1\par\vskip.5\baselineskip}
\renewcommand{\algorithmicrequire}{\textbf{Input:}}
\renewcommand{\algorithmicensure}{\textbf{Output:}}
\begin{algorithm}[t]
  \caption{ATDS algorithm with SLRA}
  \label{alg:ovralg}
  \begin{algorithmic}
    \Require $\mbA \in \btK \Cnn$ \textcolor{black}{(Set $\varepsilon > 0$ and $t_{\mathrm{max}} \in \N$ for stopping condition.)}
    \Ensure a common exact diagonalizer $\tS \in \Cnn$ of $(\mbA \approx)\tbA \in \SD$
    \State $\hA(0) = \Xi(\mbA)$; $t = 0$
    \While{$\norm{\hA(t) - P_{\mfL_{n^2 - n}}(\hA(t))}_F > \varepsilon$ \textcolor{black}{and $t < t_{\mathrm{max}}$}}
    \State $\hA(t + 1) = P_{\Xi(\btK \Cnn)} \circ P_{\mfL_{n^2 - n}}(\hA(t))$
    \State $t = t + 1$
    \EndWhile
    \textcolor{black}{
      \State $t_{\mathrm{end}} \coloneqq t$
      \If{$P_{\Xi^{-1}(\hA(t_{\mathrm{end}}))}(\mbA) \in \SD$}
      \Comment{can be judged with Fact~\ref{fct:ext}.}
      \State $\tbA = P_{\Xi^{-1}(\hA(t_{\mathrm{end}}))}(\mbA)$
      \State $\tS =$ \Call{DODO}{$\tbA, n$}
      \Comment{see Algorithm~\ref{alg:dodo} for \Call{DODO}{}.}
      \Else 
      \Comment{see Remark~\ref{rmk:prpalg}(e).}
      \State $\tS =$ \Call{PCD}{$P_{\Xi^{-1}(\hA(t_{\mathrm{end}}))}(\mbA)$}
      \Comment{see below for \Call{PCD}{}.}
      \State $\tbA = P_{\SD(\tS)}(\mbA)$.
      \EndIf
      \algrule
      \State Assumption~1: at least one of $X_{k} \in \Cnn\ (k = 1, \ldots , K)$ is diagonalizable.
      \Function{PCD}{$X_{1}, \ldots , X_{K}$}
      \Comment{PCD: "Pseudo Common Diagonalizer"}
      \State $\mathfrak{I} \coloneqq \{k \in \{1, \ldots , K\} \mid X_{k} \text{ is diagonalizable.}\}$
      \State Compute $\tS_{k} \in \Cnn\ (k \in \mathfrak{I})$ s.t. $\tS_{k}^{-1} X_{k} \tS_{k}$ is diagonal.
      \State Choose $l_{\ast} \in \argmin_{l \in \mathfrak{I}} \norm{\mbA - P_{\SD(\tS_{l})}(\mbA)}_{F}$.
      \State \Return $\tS_{l_{\ast}} \in \Cnn$
      \EndFunction
    }
  \end{algorithmic}
\end{algorithm}

\begin{rmk}[On Algorithm~\ref{alg:ovralg}]
  \label{rmk:prpalg}
  \ 
  \begin{enumerate}[label=(\alph*)]
    \item \textcolor{black}{
        The initial guess $\hA(0) = \Xi(\mbA)$ can be computed without
        multiplications because each $n$-by-$n$ block $\check{A}_{i,j}^{(k)}
        \in \Cnn\ (i, j = 1, \ldots , n)$ of $I_{n} \otimes A_{k} -
        A_{k}^{\top} \otimes I_{n}\ (k = 1, \ldots , K)$ is given by
        $\check{A}_{i,j}^{(k)} = A_{k} - a_{i,i}^{(k)} I_{n}$ if $i = j$;
        otherwise $\check{A}_{i,j}^{(k)} = -a_{j,i}^{(k)} I_{n}$, where $A_{k}
        \coloneqq [a^{(k)}_{i,j}]_{i,j=1}^{n}$.
      }
    \item The projection $P_{\mfL_{n^2 - n}}(\hX)$ in~\eqref{eq:prjs} can be
      computed with the truncated Singular Value Decomposition~(SVD) of $\hX$~(see
      the Schmidt approximation theorem in,
      e.g.,~\cite[Theorem~3]{Ben-Israel2003}).
      \textcolor{black}{
        The SVD of $\hA(t) (\in \C^{Kn^{2}\times
        n^{2}})\ (t = 0, 1, \ldots)$ is certainly dominant in the computation time of Algorithm~\ref{alg:ovralg} although efficient SVD algorithms for a large
        matrix have been studied extensively~(see~\cite{Dongarra2018} and the
        references therein).
      }
      We also remark that $P_{\mfL_{n^2 - n}}(\hX)$ in~\eqref{eq:prjs} is determined
      uniquely except in a very special case where the $(n^2 - n)$th and $(n^2
      - n + 1)$st singular values of $\hX$ happen to
      coincide~\cite[Theorem~3]{Ben-Israel2003}.
    \item The projection $P_{\Xi(\btK \Cnn)}([\hX^\top_1, \ldots ,
      \hX^\top_K]^\top) \eqqcolon [\hZ^\top_1, \ldots , \hZ^\top_K]^\top\
      (\hX_k, \hZ_k \in \C^{n^2 \times n^2},\ k = 1,\ldots , K)$
      in~\eqref{eq:prjs} can be computed by assigning, to $\hZ_{k}$, the orthogonal projection
      of $\hX_k$ onto $\mathcal{M} \coloneqq \spn\{\hE(i, j) \coloneqq I_n
        \otimes E(i,j)-E^\top(i,j)\otimes
      I_n\mid (i,j)\in\{1,\ldots ,n\}^2 \setminus \{(n,n)\}\}$,
      where $E(i,j) = [e^{(i, j)}_{p, q}]\in\Cnn$ is given by $e^{(i, j)}_{p,q}
      = 1$ if $(p, q) = (i, j)$; otherwise $0$.
      \textcolor{black}{
        Moreover, the orthogonal projections $P_{\mathcal{M}}(\hX_{k})$ can be
        computed efficiently by using the sparsity of $\hE(i, j)\ (i, j = 1,
        \ldots , n)$ and $P_{\mathcal{M}}(\hX_{k}) = P_{\spn\{\hE(i, i) \mid i
        \neq n\}}(\hX_{k}) + \sum_{1 \le s \neq t \le n} P_{\spn\{\hE(s,
        t)\}}(\hX_{k})$, where $P_{\spn\{\hE(i, i) \mid i \neq n\}}(\hX_{k})$
        and $P_{\spn\{\hE(s, t)\}}(\hX_{k})\ (s, t = 1, \ldots , n; s \neq t)$
        are the orthogonal projections onto $\spn\{\hE(i, i) \mid i \neq n\}$
        and $\spn\{\hE(s, t)\}$, respectively.
      }
    \item Alternating projection algorithm
      \textcolor{black}{
        used in Step~1'a
      }
      is a powerful tool to solve feasibility problems. Even for nonconvex
      feasibility problems, the algorithm has a guarantee to converge
      locally~\cite{Lewis2009,Noll2016} to a point in the intersection and has
      been used extensively for finding a point, near the initial guess, in the
      intersection, e.g., phase retrieval~\cite{Bauschke2002}.
    \item \textcolor{black}{
        Suppose $S \in \Cnn$ is invertible. Since $\SD(S) \coloneqq \{\mbX \in
        \SD \mid \mbX \text{ has a common diagonalizer } S\}$ is a subspace of
        $\btK \Cnn$, the orthogonal projection of $\mbA$ onto $\SD(S)$, i.e.,
        $P_{\SD(S)}(\mbA) \coloneqq (Z_{1}, \ldots , Z_{K}) \in \SD(S)$, is
        well-defined and can be computed as\footnotemark[7] 
        \begin{align*}
          \hspace{-6mm}
          Z_{k} \coloneqq \vec^{-1}(P_{\range(S^{-\top}\odot S)}(\vec(A_{k})))\ (k = 1, \ldots , K),
        \end{align*}
        where $P_{\range(S^{-\top}\odot S)}(\vec(A_{k}))$ are the orthogonal
        projections onto $\range(S^{-\top}\odot S) \subset \C^{n^{2}}$.  Even for
        the cases of $P_{\Xi^{-1}(\hA(t_{\mathrm{end}}))}(\mbA) \not \in \SD$
        in Algorithm~2, we can obtain $\tbA = P_{\SD(\tS)}(\mbA) \in \SD$, via the
        function~\texttt{PCD}, if at least one matrix in the tuple $P_{\Xi^{-1}(\hA(t_{\mathrm{end}}))}(\mbA)$ is
        diagonalizable.
      }
      \footnotetext[7]{
        \textcolor{black}{
          This is verified by the well-known identity $\vec(A \diag(\mbd) C) =
          (C^{\top} \odot A) \mbd\ (A, C \in \Cnn; \mbd \in \Cn)$~(see,
          e.g.,~\cite{Ma2010}) which is also used in the proof
          in~\ref{prf:XiAndSD}.
        }
      }
  \end{enumerate}
\end{rmk}

%% file: Sections/nmex.tex
\section{Numerical Experiments}
\label{sec:nuex}
\begin{figure*}[t]
  \begin{tabular}{ccc}
    \begin{minipage}[t]{0.33\linewidth}
      \centering
      \centerline{\includegraphics[width=55mm]{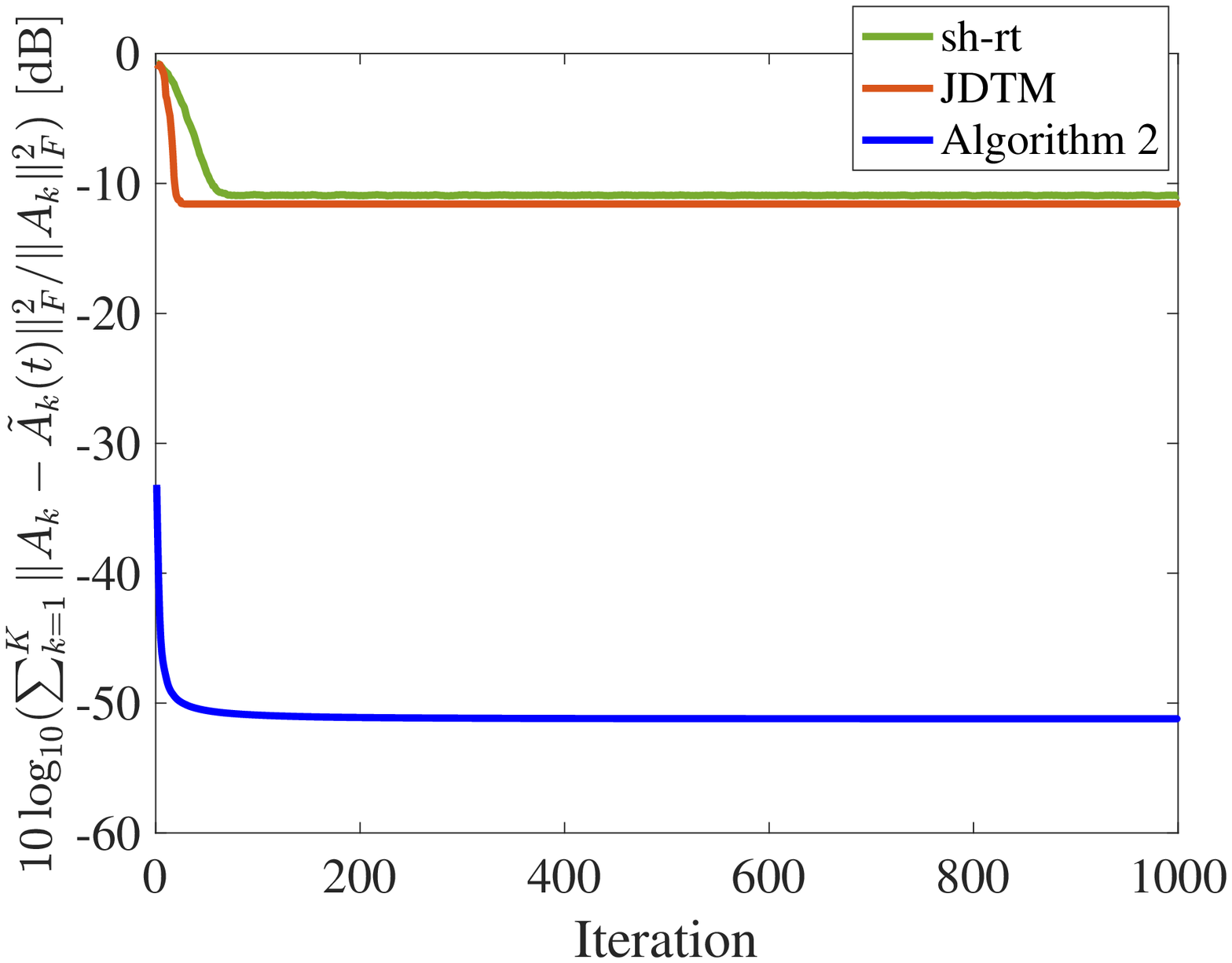}}
    \end{minipage}
    \begin{minipage}[t]{0.33\linewidth}
      \centering
      \centerline{\includegraphics[width=55mm]{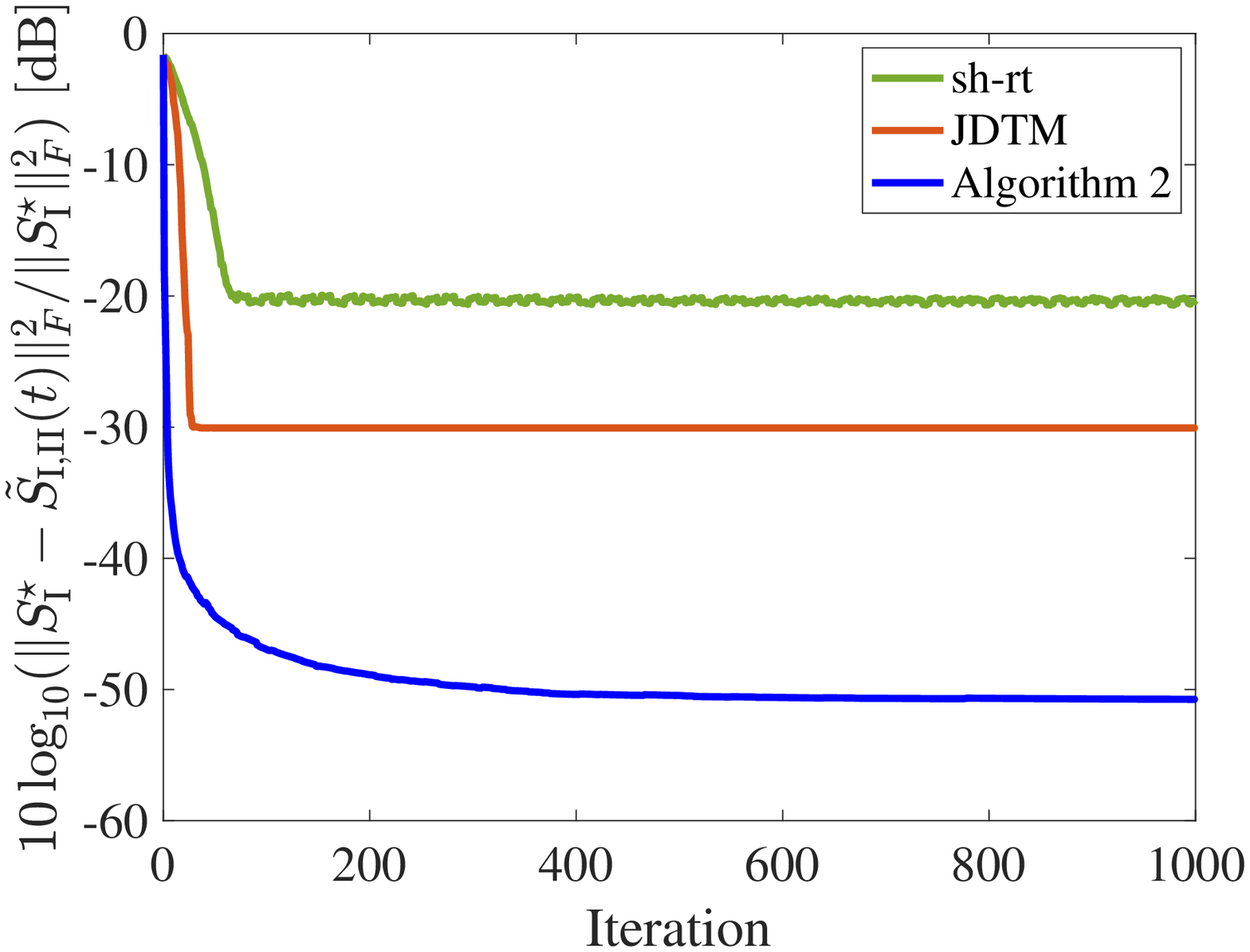}}
    \end{minipage}
    \begin{minipage}[t]{0.33\linewidth}
      \centering
      \centerline{\includegraphics[width=55mm]{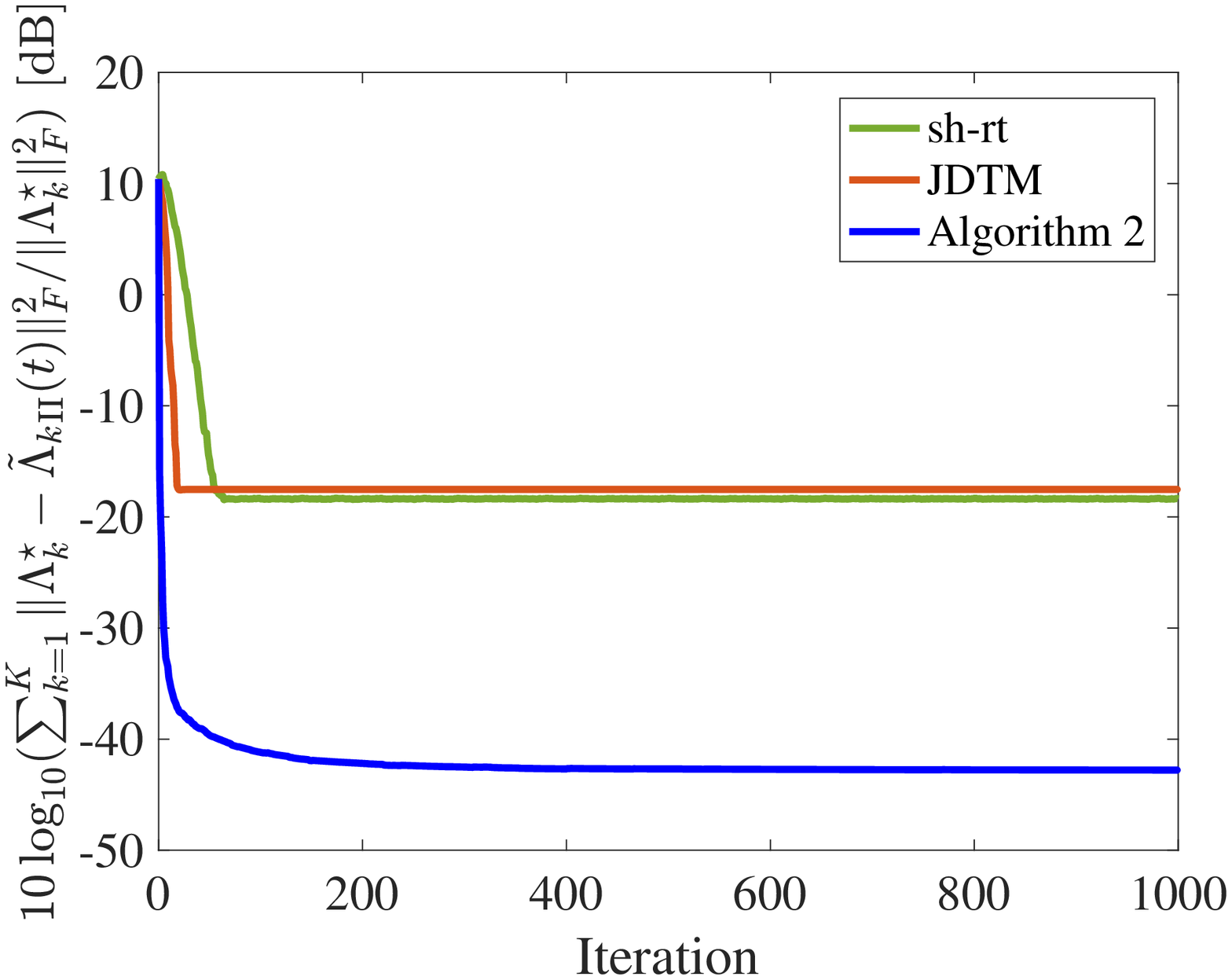}}
    \end{minipage}
  \end{tabular}
  \caption{
    \textcolor{black}{
      Transition of approximation errors in the case of SNR $50$ dB and $\kappa
      = 50$.
    }
  }
  \label{fig:iter}
\end{figure*}

\begin{table*}[t]
  \centering
  \caption{
    \textcolor{black}{
      Computation time, number of iterations, until
      $\norm{S^{\star}_{\mathrm{I}}-\tS_{\mathrm{I,II}}(t)}_{F}^{2}/\norm{S_{\mathrm{I}}^{\star}}_{F}^{2} \le 10^{-2},
      10^{-3}$, or $10^{-4}$, and rate of successful trials (from left to
      right) in the case of SNR $50$ dB and $\kappa = 50$. Note that we
      use "-" if the rate of successful trials is smaller than $10 \%$.
    }
  }
  \label{tbl:cmptime2}
  \begin{tabular}{llllllllll}
    \hline
    $\norm{S^{\star}_{\mathrm{I}}-\tS_{\mathrm{I,II}}(t)}_{F}^{2}/\norm{S_{\mathrm{I}}^{\star}}_{F}^{2}$ & \multicolumn{3}{l}{$\le 10^{-2}$} & \multicolumn{3}{l}{$\le 10^{-3}$} & \multicolumn{3}{l}{$\le 10^{-4}$} \\
    sh-rt~\cite{Fu2006} & $4.298$ ms & $52.69$ & $75 \%$ & - & - & $8 \%$ & - & - & $2 \%$ \\
    JDTM~\cite{Luciani2014} & $\mathbf{1.552}$ {\bf ms} & $18.85$ & $100 \%$ & $\mathbf{1.873}$ {\bf ms} & $22.99$ & $72 \%$ & - & - & $8 \%$ \\
    Algorithm~\ref{alg:ovralg} & $4.004$ ms & $\mathbf{1.610}$ & $100 \%$ & $5.212$ ms & $\mathbf{3.540}$ & $\mathbf{100 \%}$ & $12.32$ ms & $\mathbf{14.79}$ & $\mathbf{100 \%}$ \\
    \hline
  \end{tabular}
\end{table*}

\begin{figure*}[t]
  \begin{tabular}{ccc}
    \begin{minipage}[t]{0.33\linewidth}
      \centering
      \centerline{\includegraphics[width=55mm]{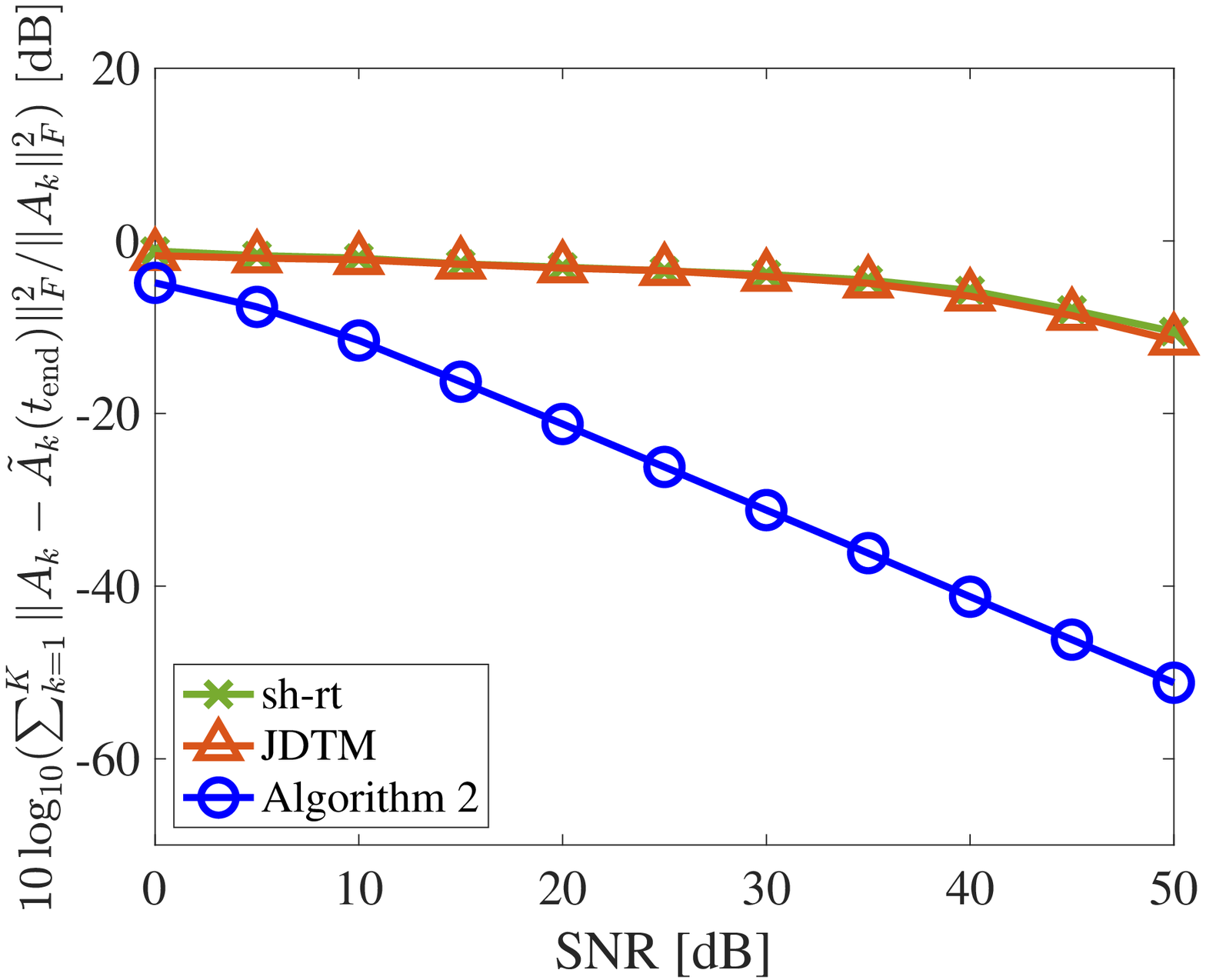}}
    \end{minipage}
    \begin{minipage}[t]{0.33\linewidth}
      \centering
      \centerline{\includegraphics[width=55mm]{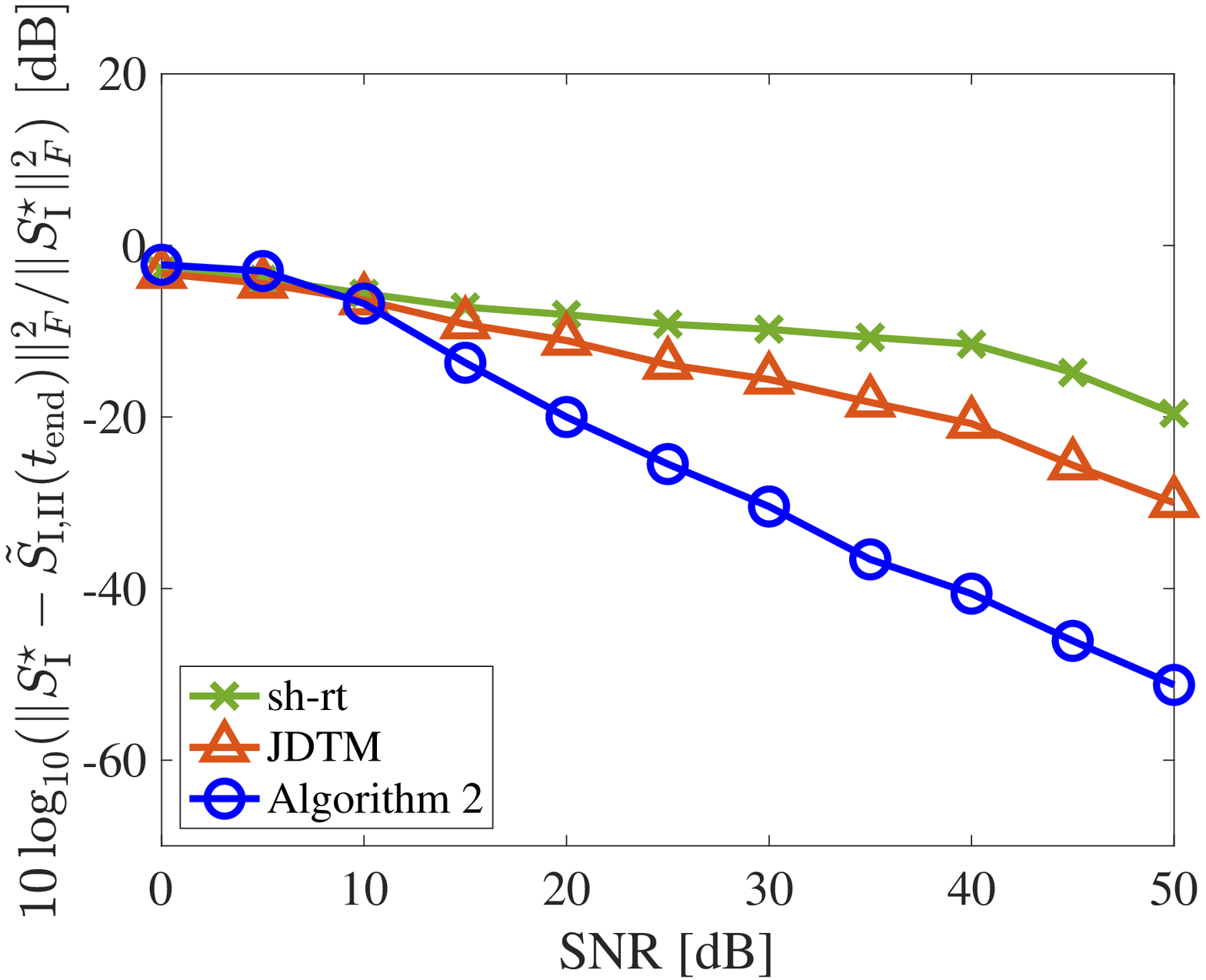}}
    \end{minipage}
    \begin{minipage}[t]{0.33\linewidth}
      \centering
      \centerline{\includegraphics[width=55mm]{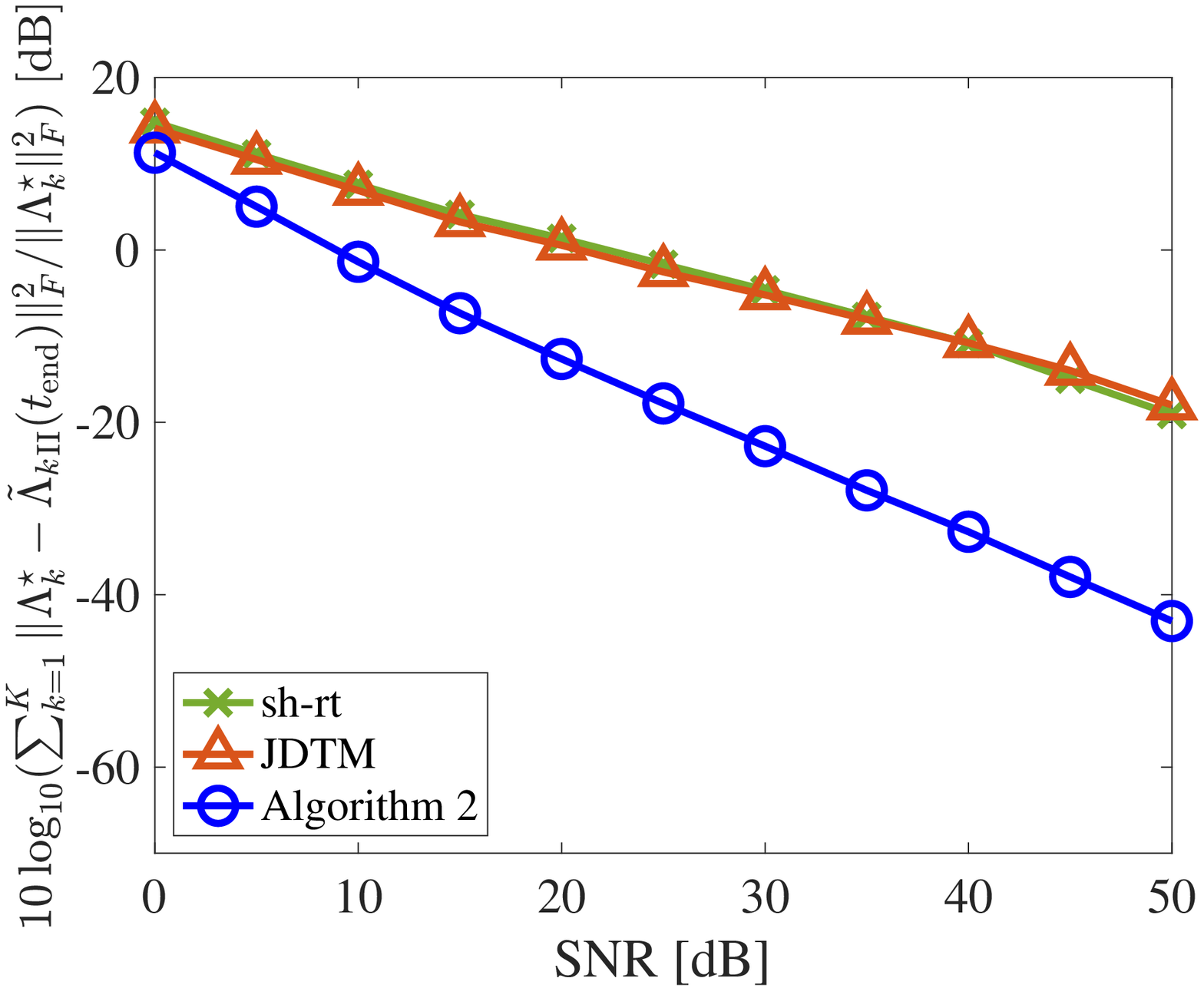}}
    \end{minipage}
  \end{tabular}
  \caption{
    Comparison of ASD algorithms for $S^{\star}$ of condition number $\kappa =
    50$.
  }
  \label{fig:asd2}
\end{figure*}

\begin{table*}[t]
  \centering
  \caption{
    \textcolor{black}{
      Computation time, number of iterations, until
      $\norm{S^{\star}_{\mathrm{I}}-\tS_{\mathrm{I,II}}(t)}_{F}^{2}/\norm{S_{\mathrm{I}}^{\star}}_{F}^{2} \le 10^{-2},
      10^{-3}$, or $10^{-4}$, and rate of successful trials (from left to
      right) in the case of SNR $50$ dB and $\kappa = 5$.
    }
  }
  \label{tbl:cmptime1}
  \begin{tabular}{llllllllll}
    \hline
    $\norm{S^{\star}_{\mathrm{I}}-\tS_{\mathrm{I,II}}(t)}_{F}^{2}/\norm{S_{\mathrm{I}}^{\star}}_{F}^{2}$ & \multicolumn{3}{l}{$\le 10^{-2}$} & \multicolumn{3}{l}{$\le 10^{-3}$} & \multicolumn{3}{l}{$\le 10^{-4}$} \\
    sh-rt~\cite{Fu2006} & $2.139$ ms & $21.27$ & $100 \%$ & $2.542$ ms & $26.83$ & $100 \%$ & $2.777$ ms & $29.72$ & $100 \%$ \\
    JDTM~\cite{Luciani2014} & $\mathbf{1.227}$ {\bf ms} & $14.09$ & $100 \%$ & $\mathbf{1.425}$ {\bf ms} & $17.13$ & $100 \%$ & $\mathbf{1.603}$ {\bf ms} & $19.06$ & $100 \%$ \\
    Algorithm~\ref{alg:ovralg} & $3.208$ ms & $\mathbf{1.000}$ & $100 \%$ & $3.334$ ms & $\mathbf{1.000}$ & $100 \%$ & $3.338$ ms & $\mathbf{1.000}$ & $100 \%$ \\
    \hline
  \end{tabular}
\end{table*}

\begin{figure*}[t]
  \begin{tabular}{ccc}
    \begin{minipage}[t]{0.33\linewidth}
      \centering
      \centerline{\includegraphics[width=55mm]{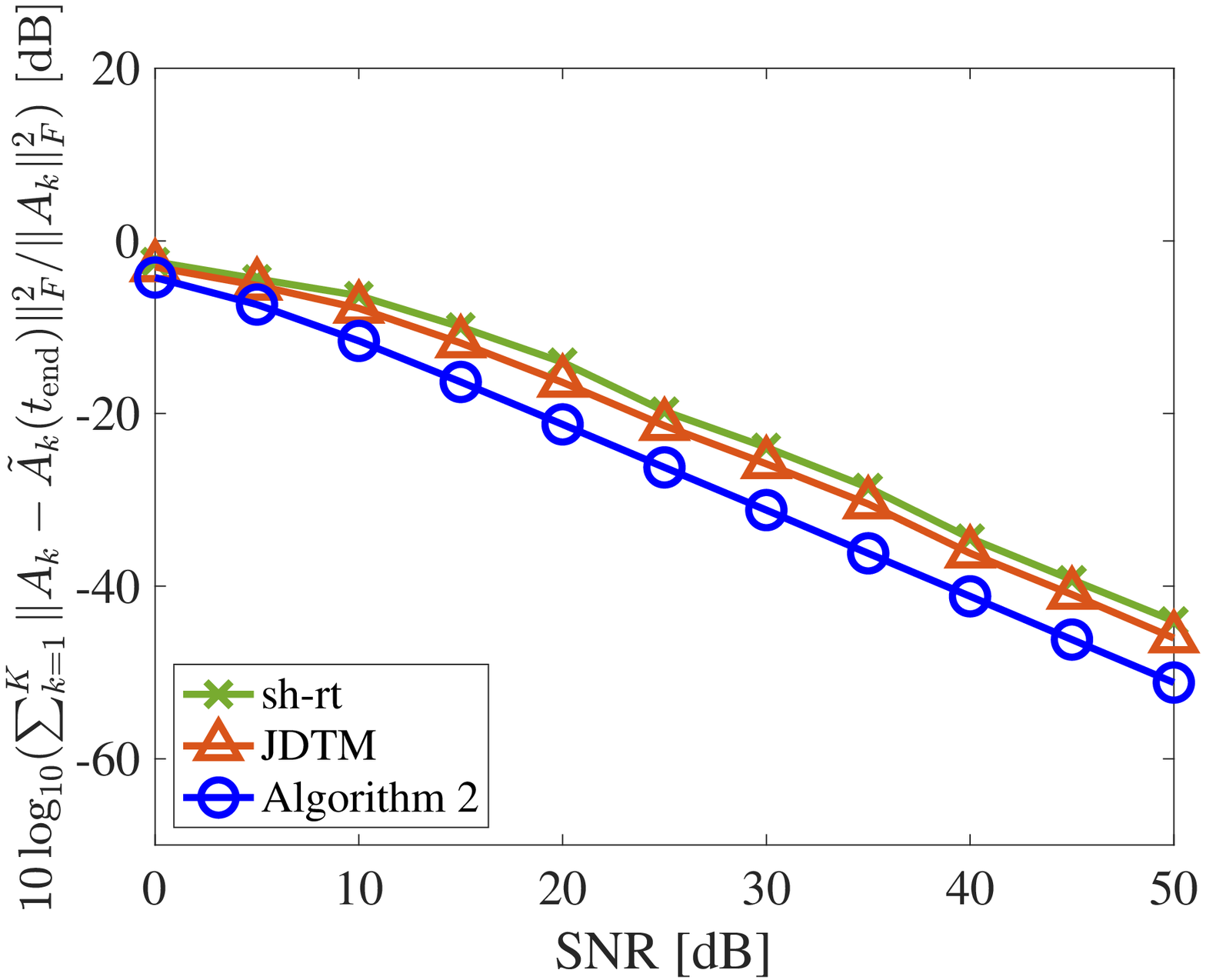}}
    \end{minipage}
    \begin{minipage}[t]{0.33\linewidth}
      \centering
      \centerline{\includegraphics[width=55mm]{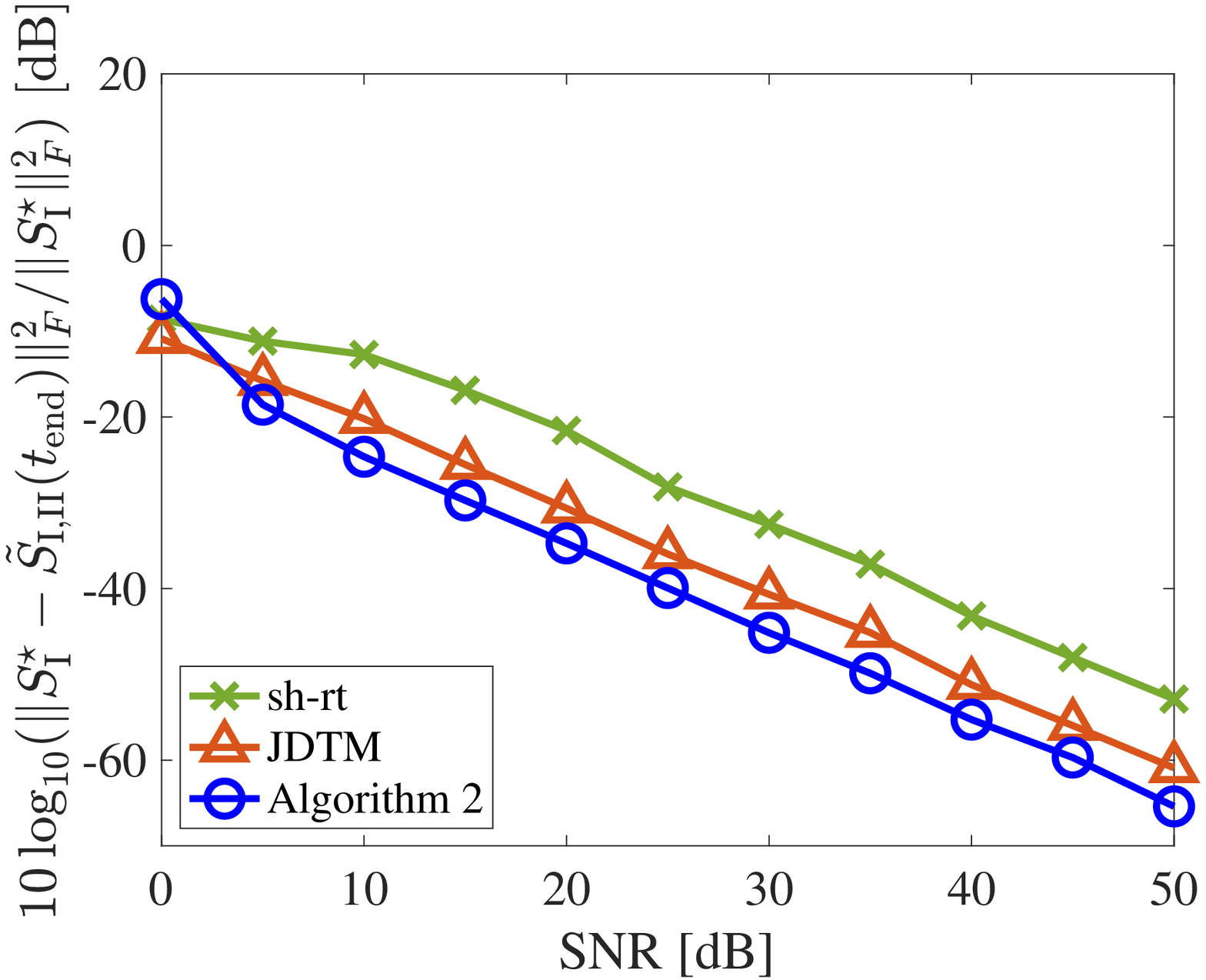}}
    \end{minipage}
    \begin{minipage}[t]{0.33\linewidth}
      \centering
      \centerline{\includegraphics[width=55mm]{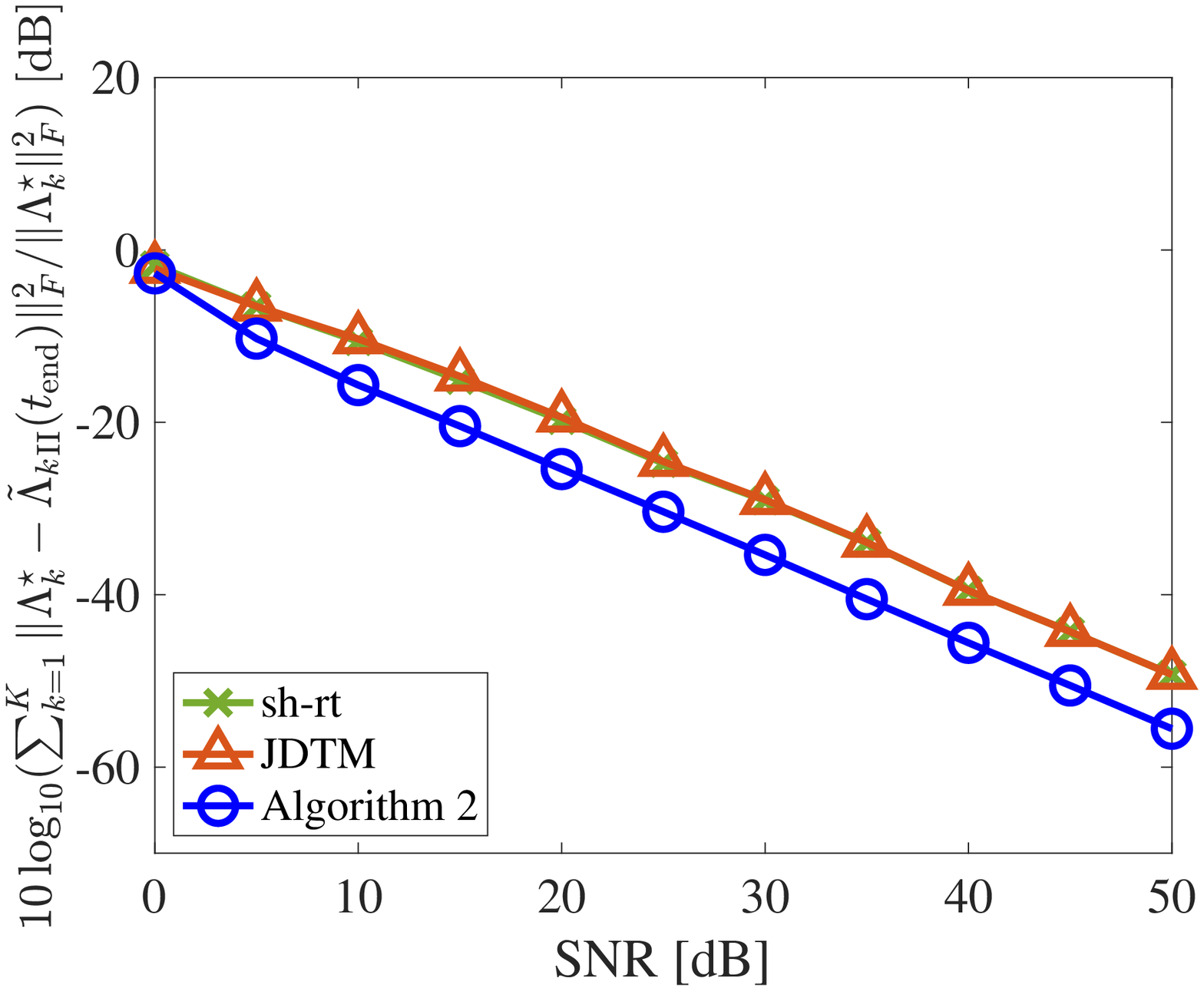}}
    \end{minipage}
  \end{tabular}
  \caption{
    Comparison of ASD algorithms for $S^{\star}$ of condition number $\kappa =
    5$.
  }
  \label{fig:asd1}
\end{figure*}

\textcolor{black}{
  To see the numerical performance of Algorithm~\ref{alg:ovralg}, in comparison
  to the two Jacobi-like methods~(sh-rt~\cite{Fu2006} and
  JDTM~\cite{Luciani2014}), under several
  conditions~(e.g., noise levels, condition numbers of an ideal common
  diagonalizer), we conduct numerical experiments for a perturbed version
  $\mbA = (A_1, \ldots , A_K) \in \btK \Rnn\ (n = 5, K = 20)$ of $\mbA^{\star} = (A^{\star}_{1}, \ldots , A^{\star}_{K}) \coloneqq (S^\star \Lambda_k^\star (S^\star)^{-1})_{k=1}^{K} \in \SD$ with $A_k \coloneqq
  A_{k}^{\star} + \sigma N_k \ (k = 1, \ldots , K)$, where the
  diagonal entries of diagonal matrices $\Lambda_k^\star \in \Rnn$ and all the entries of $N_k \in \Rnn$ are drawn
  from the standard normal distribution $\mcN(0,1)$ and $\sigma > 0$ is used to define the Signal
  to Noise Ratio~(SNR).  To conduct numerical
  experiments for $S^{\star} \in \Rnn$ of various condition numbers, say
  $\kappa > 1$, we design $S^{\star}$ by replacing singular values, of a matrix
  whose entries are drawn from the standard normal distribution $\mcN(0,1)$, with
  $\sigma_{i}(S^{\star}) \coloneqq (\kappa - 1) (n - i) / (n - 1) + 1 (i = 1,
  \ldots , n)$ implying thus $\sigma_{1}(S^{\star}) / \sigma_{n}(S^{\star}) =
  \kappa$.
  For Algorithm~\ref{alg:ovralg}, we use, as the $t$th estimates, $\tS(t) \coloneqq $ \texttt{DODO}($P_{\Xi^{-1}(\hA(t))}(\mbA), n$) if $P_{\Xi^{-1}(\hA(t))}(\mbA) \in \SD$; otherwise $\tS(t) \coloneqq$ \texttt{PCD}($P_{\Xi^{-1}(\hA(t))}(\mbA)$)~(see function \texttt{PCD} in Algorithm~\ref{alg:ovralg}),\footnotemark[8] $\tbA(t) = (\tA_{1}(t), \ldots , \tA_{K}(t)) \coloneqq P_{\Xi^{-1}(\hA(t))}(\mbA)$ if $P_{\Xi^{-1}(\hA(t))}(\mbA) \in \SD$; otherwise $\tbA(t) \coloneqq P_{\SD(\tS(t))}(\mbA)$, and $\tLam_{k}(t) \coloneqq (\tS(t))^{-1} \tA_{k}(t) \tS(t) \ (k = 1, \ldots , K; t = 0, 1, \ldots)$.
  For the
  Jacobi-like methods, on the other hand, we use, as the $t$th estimates, $\tS(t) \coloneqq \breve{S}(t)$, $\tLam_{k}(t) \coloneqq
  \breve{\Lambda}_{k}(t)$, where $\breve{S}(t)$ and $\breve{\Lambda}_{k}(t)$ are respectively the $t$th updates of $\breve{S}$ and $\breve{\Lambda}$ in the Jacobi-like methods~(sh-rt and JDTM; see Section~\ref{ssec:jcb} for sh-rt), and $\tA_{k}(t)
  \coloneqq \tS(t) \tLam_{k}(t) (\tS(t))^{-1}\ (k = 1, \ldots , K)$.
  The Jacobi-like
  methods are terminated when the number of iteration
  exceeds $2 \times 10^{4}$ or when
  $\abs{f_\mbA(\tS(t))-f_{\mbA}(\tS(t-1))}/\abs{f_{\mbA}(\tS(t))} \le 10^{-6}\
  (t \in \N)$, where $f_\mbA(\tS(t)) = \sum_{k = 1}^K
  \off((\tS(t))^{-1} A_k \tS(t))$. We choose $\varepsilon = 10^{-6}$ and
  $t_{\mathrm{max}} = 2 \times 10^{4}$ in Algorithm~\ref{alg:ovralg}. For each algorithm, we use $t_{\mathrm{end}} \in \N$ to indicate the iteration when the algorithm is terminated. We
  evaluate the approximation errors of $\tS(t)$, $(\tLam_{1}(t), \ldots , \tLam_{K}(t))$, $\tbA(t) = (\tA_{1}(t), \ldots , \tA_{k}(t))$ by 
  $\norm{S^{\star}_{\mathrm{I}}-\tS_{\mathrm{I, II}}(t)}^{2}_{F} / \norm{S^{\star}_{\mathrm{I}}}^{2}_{F}$,
  $\sum_{k=1}^{K} \norm{\Lambda^{\star}_{k}-\mbox{$\tLam_{k}$}_{\mathrm{II}}(t)}^{2}_{F} /
  \norm{\Lambda^{\star}_{k}}^{2}_{F}$, and $\sum_{k=1}^{K}
  \norm{A_{k}-\tA_{k}(t)}^{2}_{F} / \norm{A_{k}}^{2}_{F}$, respectively, where (i)~$S^{\star}_{\mathrm{I}}$ stands for the column-wise normalization of $S^{\star}$ and $\tS_{\mathrm{I, II}}(t)$ stands for the column-wise permutation applied to achieve the best approximation to $S_{I}^{\star}$ after the column-wise normalization of $\tS(t)$, and (ii)~$\mbox{$\tLam_{k}$}_{\mathrm{II}}(t)$ stands for the diagonal matrix after applying the corresponding permutation for $\tS_{\mathrm{I,II}}(t)$ to diagonal entries of $\tLam_{k}(t)$.   
}

\footnotetext[8]{
  \textcolor{black}{
    In our numerical experiments, we have not seen any exceptional case where
    Assumption~1 for $P_{\Xi^{-1}(\hA(t))}(\mbA)$ is not satisfied.
  }
}

\textcolor{black}{
  We conducted numerical experiments on Intel Core i7-8559U running at $2.7$ GHz
  with $4$ cores and $32$ GB of main memory. By using Matlab, we implemented all the ASD
  algorithms by ourselves.  We measured the computation times of
  all the ASD algorithms by Matlab's \texttt{tic/toc} functions. 
  We compared the computation time and the number of iterations, of all the algorithms, until
  $\norm{S^{\star}_{\mathrm{I}}-\tS_{\mathrm{I,II}}(t)}_{F}^{2}/\norm{S^{\star}_{\mathrm{I}}}_{F}^{2} \le 10^{-2},
  10^{-3}$, and $10^{-4}\ (t = 0, 1, \ldots)$ as shown in Table~2 and~3.
  For each ASD algorithm, the successful trial means a trial where the
  algorithm succeeds in achieving smaller approximation error of $\tS(t)$ than prescribed values $10^{-2}, 10^{-3}$, and $10^{-4}$ before termination of each algorithm.
}

\textcolor{black}{
  Since it is reported in~\cite{Luciani2014} that the Jacobi-like methods tend to suffer from the cases where the input tuple is a perturbed version of a simultaneously diagonalizable tuple with a common diagonalizer of large condition number, we first compared all the ASD algorithms for $\kappa = 50$. Figure~\ref{fig:iter} depicts the transition of the mean values, of the relative squared errors of $\tbA(t), \tS(t)$ and $(\tLam_{1}(t), \ldots , \tLam_{K}(t))$ after proper column-wise normalization/permutation for $\tS(t)$ and $(\tLam_{k}(t))_{k=1}^{K}$, over $100$ trials in the case of SNR $50$ dB.
  Figure~1 illustrates that, compared with the Jacobi-like methods, Algorithm~\ref{alg:ovralg} achieves estimations (i)~$\tbA(t) \in \SD$, of $\mbA^{\star}$, closer to $\mbA$, (ii)~$\tS_{\mathrm{I,II}}(t)$, of $S_{\mathrm{I}}^{\star}$, closer to $S_{\mathrm{I}}^{\star}$, and (iii)~$(\mbox{$\tLam_{k}$}_{\mathrm{II}}(t))_{k=1}^{K}$, of $(\Lambda^{\star}_{k})_{k=1}^{K}$, closer to $(\Lambda^{\star}_{k})_{k=1}^{K}$ with smaller number of iterations.
  Table~\ref{tbl:cmptime2} depicts the mean values of (i)~the computation times
  and (ii)~the numbers of iterations taken until
  $\norm{S^{\star}_{\mathrm{I}}-\tS_{\mathrm{I,II}}(t)}_{F}^{2}/\norm{S_{\mathrm{I}}^{\star}}_{F}^{2} \le 10^{-2},
  10^{-3}$, or $10^{-4}$ over successful trials in $100$ trials and (iii)~the rates
  of successful trials over $100$ trials.
  This result shows that Algorithm~\ref{alg:ovralg} takes around
  $3$ times longer computation time than JDTM but its $\tS(t_{\mathrm{end}})$ achieves
  the prescribed conditions even for the trials where the Jacobi-like
  methods fail to achieve the prescribed conditions. Figure~\ref{fig:asd2} depicts the mean values, of the relative squared errors of $\tbA(t_{\mathrm{end}}), \tS(t_{\mathrm{end}})$, and $(\tLam_{k}(t_{\mathrm{end}}))_{k=1}^{K}$ after proper column-wise normalization/permutation for $\tS(t)$ and $(\tLam_{k}(t))_{k=1}^{K}$, over $100$ trials in the cases of SNR from $0$ dB to $50$ dB. Figure~2 illustrates Algorithm~2 outperforms the
  Jacobi-like methods in the sense of achieving approximation errors at $t_{\mathrm{end}}$ especially when SNR is higher than $10$ dB.
}

\textcolor{black}{
  We also compared all the algorithms in the case of $\kappa = 5$.
  All the values in Table~\ref{tbl:cmptime1} and Figure~\ref{fig:asd1} are
  calculated by the same way as done for $\kappa = 50$.
  Figure~\ref{fig:asd1} and Table~\ref{tbl:cmptime1} show that Algorithm~2 takes around $3$ times longer computation time than JDTM but can outperform the
  Jacobi-like methods in the sense of achieving approximation errors at $t_{\mathrm{end}}$ in SNR from $5$ dB to $50$ dB.
  From these experiments, we see that Algorithm~\ref{alg:ovralg} is robust against wider range of condition numbers of $S^{\star}$ than the Jacobi-like methods.
}

%% file: Sections/cncl.tex
\section{\textcolor{black}{Concluding Remarks}}
\label{sec:cncl}
\textcolor{black}{
In this paper, for the approximate simultaneous diagonalization of $\mbA \in \btK \Cnn$, we
newly presented the Approximate-Then-Diagonalize-Simultaneously~(ATDS)
algorithm by solving a certain Structured Low-Rank Approximation~(SLRA). For $\mbA \in \SD$, the
proposed ATDS algorithm has a guarantee to find its
common exact diagonalizer unlike the Jacobi-like methods.
Numerical experiments show that, at the expense of reasonable computational time, the
proposed ATDS algorithm achieves better approximations to the desired information in Problem~2 than the Jacobi-like methods for almost all SNR as well as both for small and large condition numbers of $S^{\star}$.
}

\textcolor{black}{
  The reduction of the computational cost for SVD of $\hA(t) \in \C^{K n^{2}\times n^{2}}\ (t = 0, 1, \ldots)$ in Algorithm~2 is under study. We are also studying applications of the proposed ATDS algorithm to certain signal processing problems.\footnote[9]{\textcolor{black}{A partial result for applications of the ATDS algorithm to Canonical
  Polyadic~(CP) tensor decomposition was presented at a
conference~\cite{Akema2020}.}} These will be reported elsewhere.
}

%% file: Sections/apdx.tex
\section{\textcolor{black}{Nonconvexity of $\mfD^{\mathrm{abl}}_{2,2}$}}
\label{apd:ncsd}
\textcolor{black}{
  The following simple example shows the nonconvexity of $\mfD^{\mathrm{abl}}_{2,2} \subsetneq \C^{2\times 2} \times \C^{2\times 2}$.
  Let $A_{1} \coloneqq 
  \big[
    \begin{smallmatrix}
      1 & 1 \\
      -1 & 1
    \end{smallmatrix}
  \big], A_{2} \coloneqq 2 A_{1}, B_{1} \coloneqq
  \big[
    \begin{smallmatrix}
      1 & 1 \\
      1 & 1
    \end{smallmatrix}
  \big]$, and $B_{2} \coloneqq 2 B_{1}$. Since $A_{1}, A_{2}, B_{1}$
  and $B_{2}$ are normal, these are respectively
  diagonalizable~\cite[Theorem~2.5.3]{Horn2013}. Moreover, by $A_{1}
  A_{2} = A_{2} A_{1}$ and $B_{1} B_{2} = B_{2} B_{1}$, we see $\mbA \coloneqq (A_{1}, A_{2}), \mbB \coloneqq (B_{1}, B_{2}) \in \mfD^{\mathrm{abl}}_{2,2}$ (see Fact~\ref{fct:ext} in Section~\ref{ssec:dfnfct}). Below, we will show that $\mbA / 2 + \mbB / 2 \not \in \mfD^{\mathrm{abl}}_{2,2}$.
  Since every eigenvector of $(A_{1} + B_{1})/2 =
  \big[
    \begin{smallmatrix}
      1 & 1 \\
      0 & 1
    \end{smallmatrix}
  \big]$ is given by $\mbv = 
  \big[
    \begin{smallmatrix}
      v \\
      0
    \end{smallmatrix}
  \big]\ (v \in \C \setminus \{0\})$, the eigenspace of $
  \big[
    \begin{smallmatrix}
      1 & 1 \\
      0 & 1
    \end{smallmatrix}
  \big]$ has dimension $1$, which implies $(A_{1} + B_{1})/2$ is not
  diagonalizable and hence $\mbA / 2 + \mbB / 2 \not \in
  \mfD^{\mathrm{abl}}_{2,2}$.
}

\section{Useful Facts on Commutativity and Diagonalizability}
\label{apd:fcdo}
\begin{fct}[On commutativity]
  \label{fct:cmm}
  \ 
  \begin{enumerate}[label=(\alph*)]
    \item \label{fct:cmmb} For a given $X \in \Cnn$, the set of all matrices
      which commute with $X$ is a subspace of $\Cnn$ with dimension at least
      $n$; the dimension is equal to $n$ if and only if each eigenvalue of $X$
      has geometric multiplicity $1$, i.e., $X$ is
      nonderogatory~\cite[Corollary~4.4.15]{Horn1991}.
    \item \label{fct:cmmc} For a given $X \in \Cnn$, $X$ is nonderogatory if
      and only if every matrix $Y \in \Cnn$ which commutes with $X$ can be
      expressed as $Y = c_0 I + c_1 X + c_2 X^2 + \cdots + c_{n - 1} X^{n - 1}$
      for some $c_0, c_1, \ldots , c_{n - 1} \in \C$, i.e., $Y$ is a polynomial
      in $X$~\cite[Corollary~4.4.18]{Horn1991}.
    \item Let $\lambda_1, \ldots , \lambda_d \in \C$ be distinct and let
      $\Lambda = \lambda_1 I_{n_1} \oplus \cdots \oplus \lambda_d I_{n_d} \in
      \Cnn$. Then, for $X = [X_{i, j}]_{i, j = 1}^d \in \Cnn\ (X_{i, j} \in
      \C^{n_i \times n_j})$, $\Lambda X = X \Lambda \Leftrightarrow X = X_{1,
      1} \oplus \cdots \oplus X_{d, d}$ (Note: (c) is verified by $\Lambda X
      = X \Lambda \Leftrightarrow (\forall i, j \in \{1, \ldots , d\})
      \lambda_i X_{i, j} = \lambda_j X_{i, j}$, and $(\forall i \neq j)
      \lambda_i X_{i, j} = \lambda_j X_{i, j} \Rightarrow (\forall i \neq j)
      X_{i, j} = O$).
  \end{enumerate}
\end{fct}

\begin{fct}[On diagonalizability of block diagonal matrix~{\cite[Lemma~1.3.10]{Horn2013}}]
  \label{fct:bldg}
  Suppose $X = \check{X}_1 \oplus \cdots \oplus \check{X}_l \in \Cnn$ for some
  $X_p \in \C^{n_p \times n_p}\ (n_1 + \cdots + n_l = n)$. Then $X \in \Dga
  \Leftrightarrow (\forall p \in \{1, \ldots , l\}) \check{X}_p \in
  \msD^\mathrm{abl}_{n_p}$.
\end{fct}

\section{Structured Low-Rank Approximation}
\label{apd:slra}
\emph{Structured Low-Rank Approximation~(SLRA)}~\cite{Markovsky2008} is
a problem, for a given matrix $A \in \Cmn$, a given integer $r \in (0, \min\{m,
n\})$, and a given affine subspace $\Omega \subset \Cmn$, to find a minimizer,
in all matrices of rank at most $r$, of $\norm{A - X}^2_F$.  SLRA has many
applications in signal processing~(see e.g.,~\cite{Markovsky2008} and
references therein).

For SLRA, Cadzow's algorithm~\cite{Cadzow1988} also known as alternating
projection algorithm has been used extensively while some methods to find its
local minimizer~\cite{Markovsky2005,Schost2016} or its global
one~\cite{Ottaviani2014} also have been proposed.  It is reported that, in
practice, Cadzow's algorithm finds a structured low-rank matrix close to
a given one~\cite{Cadzow1988}.

\section{Proof of Theorem~\ref{thr:XiAndSD}}
\label{prf:XiAndSD}
\ 
\begin{enumerate}[label=(\alph*)]
  \item This follows from the expression of the condition $X_k X_l-X_lX_k=O\
    (k,l=1,\ldots ,K)$ in vector form.
  \item (Proof of "$\Rightarrow$")\ Let $X_k = S \Lambda_k S^{-1}$ be such that
    $\Lambda_k \in \Cnn\ (k = 1, \ldots , K)$ are diagonal.  By using
    identities: $(A\otimes B)(C\otimes D)=AC\otimes BD\
    (A,B,C,D\in\Cnn)$~\cite[Lemma 4.2.10]{Horn1991} and $(A\otimes
    B)^{-1}=A^{-1}\otimes B^{-1}\ (A,B\in\Cnn)$~\cite[Corollary
    4.2.11]{Horn1991}, we get $I_{n} \otimes X_k - X_k^\top \otimes I_{n}
    = (S^{-\top} \otimes S)(I_{n} \otimes \Lambda_k - \Lambda_k \otimes
    I_{n})(S^{-\top}\otimes S)^{-1}$. Since $(I_{n} \otimes \Lambda_k -
    \Lambda_k \otimes I_{n})$ is diagonal and the $((j - 1) n +
    j)$th diagonal entry is $0$ for every $j = 1, \ldots , n$, we see
    $(\forall k \in \{1, \ldots , K\})\ \nullspace(I_{n} \otimes X_k - X_k^\top
    \otimes I_{n}) \supset \range(S^{-\top} \odot S) \Leftrightarrow
    \nullspace(\Xi(\mbX)) \supset \range(S^{-\top} \odot S)$.

    (Proof of "$\Leftarrow$")\ \textcolor{black}{
      Note first that any $\mby \in \range(S^{-\top} \odot S)$ can be expressed
      with $\boldsymbol{\lambda} \coloneqq [\lambda_{i}] \in \Cn$ as $\mby =
      (S^{-\top} \odot S)\boldsymbol{\lambda}$ and that $Y \coloneqq
      \vec^{-1}(\mby) = S \diag(\lambda_{1}, \ldots , \lambda_{n}) S^{-1}$,
      where we used the well-known identity $\vec(A \diag(\mbd) C) = (C^{\top}
      \odot A) \mbd\ (A, C \in \Cnn; \mbd \in \Cn)$~(see, e.g.,~\cite[Property~1]{Ma2010}).
      By $\mby \in \nullspace(\Xi(\mbX))$ and~(a), we see $X_k Y = Y X_k\ (k = 1,
      \ldots , K)$. Choose specially $\boldsymbol{\lambda}$ such that all
      entries $\lambda_{i}\ (i = 1, \ldots , n)$ are distinct in order to
      ensure $Y$ is nonderogatory. In this case,
      Fact~\ref{fct:cmm}\ref{fct:cmmc} ensures that $X_k\ (k = 1, \ldots , K)$
      are polynomials in $Y$, say $X_{k} = \mathcal{P}_{k}(Y)$, and therefore,
      $X_{k} = S \diag(\mathcal{P}_{k}(\lambda_{1}), \ldots , \mathcal{P}_{k}(\lambda_{n}))
      S^{-1}$, which implies $\mbX \in \SD$.}
  \item Let $X_k = S \Lambda_k S^{-1}$ be such that $\Lambda_k \in \Cnn\ (k =
    1, \ldots , K)$ are diagonal. Since $\dim(\range(S^{-\top} \odot S)) = n$,
    we see from~(b) that $\dim(\nullspace(\Xi(\mbX))) \ge n \Leftrightarrow
    \rank(\Xi(\mbX)) \le n^2 - n$. Moreover, since $\Xi(\mbX) = (I_K \otimes
    (S^{-\top} \otimes S)) \Xi(\Lambda_1, \Lambda_2, \ldots , \Lambda_K)
    (S^{-\top} \otimes S)^{-1}$, we see $\rank(\Xi(\mbX)) =
    \rank(\Xi(\Lambda_1, \ldots , \Lambda_K))$.  Fact~\ref{fct:essUnq}
    ensures the simultaneous diagonalization of $\mbX$ is essentially unique
    if and only if $\Xi(\boldsymbol{\Lambda})$ has just $n$ zero vectors as
    the column vectors, i.e., $\rank(\Xi(\Lambda_1, \ldots , \Lambda_K)) =
    n^2 - n$.
  \item (Proof of "$\Rightarrow$") It follows from Fact~\ref{fct:essUnq} and~(c).
    
    (Proof of "$\Leftarrow$") Suppose that $X_1$ has $n$ distinct
    eigenvalues. Therefore, there exist $S \in \Cnn$ and a diagonal matrix
    $\Lambda_1 \in \Cnn$ such that $X_1 = S \Lambda_1 S^{-1}$. By
    Fact~\ref{fct:cmm}\ref{fct:cmmb}, $\dim(\nullspace(I_n\otimes
    X_1-X_1^\top\otimes I_n))=n$ and hence $\nullspace(\Xi(\mbX))=\nullspace(I_n\otimes
    X_1-X_1^\top\otimes I_n)$. For any $k = 2, \ldots , K$, moreover, $\nullspace(I_n\otimes
    X_k-X_k^\top\otimes I_n)\supseteq\nullspace(\Xi(\mbX))=\nullspace(I_n\otimes
    X_1-X_1^\top\otimes I_n)\ni\vec(X_1)$, which ensures the commutativity of
    $X_1$ and $X_k$. By using Fact~\ref{fct:cmm}\ref{fct:cmmc}, we see that
    each $X_k$ is a certain polynomial in $X_1$. Therefore, $S^{-1}X_kS\
    (k=2,\ldots ,K)$ are diagonal.
  \item From~(d), it is sufficient to show $\rank(\Xi(\mbY)) = n^2 - n$.  Since
    $Y_l$ has $n$ distinct eigenvalues, Fact~\ref{fct:cmm}\ref{fct:cmmb} and
    $\Xi(\mbY) \in \mfL_{n^2 - n}$ ensure $n^2 - n = \rank(I_n \otimes Y_l
    - Y_l^\top \otimes I_n) \le \rank(\Xi(\mbY)) \le n^2 - n$. \qed
\end{enumerate}

\section{Proof of Proposition~\ref{prp:pjXi}}
\label{prf:pjXi}
\begin{enumerate}[label=(\alph*)]
  \item It is clear that $\mbX \in (\spn\{I_n\})^K \Leftrightarrow \Xi(\mbX) =
    O$. \\ (proof of "$\supset$") Since, for any $\mbX\in\mbX^\diamond +
    (\spn\{I_n\})^K$, $\Xi(\mbX)=\Xi(\mbX^\diamond)=\hX$ holds,
    we see $\mbX\in\Xi^{-1}(\hX)$.\\ (proof of "$\subset$") For any
    $\mbX\in\Xi^{-1}(\hX)$, we have
    $\Xi(\mbX-\mbX^\diamond)=O\Leftrightarrow\mbX-\mbX^\diamond\in
    (\spn\{I_n\})^K$, where $\mbX^\diamond\in\Xi^{-1}(\hX)$.  Therefore
    $\mbX\in\mbX^\diamond+ (\spn\{I_n\})^K$.
  \item This follows from~(a) and equalities
    \begin{align*}
      & \min_{(Y_1,\ldots ,Y_K)\in\Xi^{-1}(\hX)} \sum_{k=1}^K\norm{X_k-Y_k}_F^2 \\
      & = \min_{\substack{\alpha_1,\ldots ,\alpha_k\in\mathbb{C}}} \sum_{k=1}^K\norm{X_k-X^\diamond_k-\alpha_kI_n}^2_F \\
      & = \sum_{k=1}^K\norm{X_k - X^\diamond_k - P_{\spn\{I_n\}}(X_k-X^\diamond_k)}^2_F.
    \end{align*}
  \item For $X = [x_{i, j}] \in \Cnn$ and $Y = [y_{i, j}] \in \Cnn$,
    \begin{align*}
      & \norm{I_n \otimes X - X^\top \otimes I_n - I_n \otimes Y - Y^\top \otimes I_n}^2_F \\
      & = \sum_{i = 1}^n \norm{X - x_{i, i} I_n - Y + y_{i, i} I_n}^2_F \\ 
      & \quad + n \sum_{\substack{p, q = 1, p \neq q}}^n \abs{x_{p, q} - y_{p, q}}^2 \\
      & = n \norm{X - Y}^2_F + n \sum_{i = 1}^n \abs{x_{i, i} - y_{i, i}}^2 \\
      & \quad - 2 \abs{\trc(X - Y)}^2 + n \sum_{\substack{p, q = 1, p \neq q}}^n \abs{x_{p, q} - y_{p, q}}^2 \\
      & = 2 n \norm{X - Y}^2_F - 2 \abs{\trc(X - Y)}^2.
    \end{align*}
    Hence $\norm{\Xi(\mbX) - \hX}^2_F = \norm{\Xi(\mbX) - \Xi(Z_1,
    \ldots , Z_K)}_F^2 = \sum_{k = 1}^K 2 n \norm{X_k - Z_k}^2_F - 2 \abs{\trc(X_k
    - Z_k)}^2$. Since~(b) ensures $\trc(X_k - Z_k) = 0\ (k = 1, \ldots , K)$, we
    get $\sum_{k = 1}^K \norm{X_k - Z_k}^2_F = \norm{\Xi(\mbX) - \hX}^2_F / 2 n$.
    \qed
\end{enumerate}

\section{Proof of Proposition~\ref{prp:mnap}}
\label{prf:mnap}
This follows from $\norm{\Xi(\tbA(t)) - P_{\mfL_{n^2 - n}} (\Xi(\tbA(t)))}_F \\
\ge \norm{\Xi(\tbA(t + 1)) - P_{\mfL_{n^2 - n}} (\Xi(\tbA(t)))}_F \\ 
\ge \norm{\Xi(\tbA(t + 1)) - P_{\mfL_{n^2 - n}} (\Xi(\tbA(t + 1)))}_F$. \qed